\documentclass[msom]{informs4}

\usepackage{amsmath,amssymb,amsfonts}
\usepackage{mathtools}
\usepackage{mathrsfs}
\usepackage{bm}
\usepackage{dsfont} 

\usepackage{natbib}
\usepackage{color}
\usepackage{subfig}
\usepackage{multirow} 
\usepackage{algorithm} 
\usepackage{algorithmic} 
\usepackage{booktabs}
\usepackage{threeparttable}
\usepackage{epstopdf} 
\bibpunct[, ]{(}{)}{,}{a}{}{,}%
\def\bibfont{\small}%
\def\bibsep{\smallskipamount}%
\def\bibhang{24pt}%
\input{latexmacro}
\def\bibfont{\footnotesize}
\makeatletter
\renewcommand*\env@matrix[1][c]{\hskip -\arraycolsep
	\let\@ifnextchar\new@ifnextchar
	\array{*\c@MaxMatrixCols #1}}
\makeatother

\TheoremsNumberedThrough     
\ECRepeatTheorems

\EquationsNumberedThrough    

\MANUSCRIPTNO{}

\begin{document}


\RUNAUTHOR{Liu and Luo}

\RUNTITLE{Dynamic Dispatching and Routing with Random Demand}

\TITLE{On-Demand Delivery from Stores: Dynamic Dispatching and Routing with Random Demand}

\ARTICLEAUTHORS{%
	\AUTHOR{Sheng Liu}\AFF{Rotman School of Management, University of Toronto \\ \EMAIL{sheng.liu@rotman.utoronto.ca}}
	
	\AUTHOR{Zhixing Luo}\AFF{School of Management and Engineering, Nanjing University \\ \EMAIL{luozx.hkphd@gmail.com}}
} 

\ABSTRACT{On-demand delivery has become increasingly popular around the world. Motivated by a large grocery chain store who offers fast on-demand delivery services, we model and solve a stochastic dynamic driver dispatching and routing problem for last-mile delivery systems where on-time performance is the main target. The system operator needs to  dispatch a set of drivers and specify their delivery routes facing random demand that arrives over a fixed number of periods. The resulting stochastic dynamic program is challenging to solve due to the curse of dimensionality. We propose a novel structured approximation framework to approximate the value function via a parametrized dispatching and routing policy. We analyze the structural properties of the approximation framework and establish its performance guarantee under large-demand scenarios. We then develop efficient exact algorithms for the approximation problem based on Benders decomposition and column generation, which deliver verifiably optimal solutions within minutes. The evaluation results on a real-world data set show that our framework outperforms the current policy of the company by 36.53\% on average in terms of delivery time. We also perform several policy experiments to understand the value of dynamic dispatching and routing with varying fleet sizes and dispatch frequencies.
}%



\maketitle

%


\section{Introduction}
E-commerce has evolved rapidly and continued pushing the boundaries of digital platforms. Customers place orders online not only for household goods and electronics, but also for more perishable and time-sensitive goods such as grocery and food. The global online grocery market accounted for \$154.96 billion in 2018 and is expected to reach \$975.16 billion by 2027 \citep{og2020}. As a front runner in this market, Amazon has been providing free two-hour grocery delivery to its prime members since 2019 \citep{amzfresh2019}. In the meanwhile, meal delivery has been growing at a similar speed. The global online meal delivery market is expected to grow from \$91.41 billion in 2018 to \$182.33 billion by 2024 \citep{food2020}. 

While platforms such as Instacart and Grubhub enable small business owners to reach a broader customer base, grocery stores and restaurants try to build their own delivery capacity to maintain a reliable delivery service. McDonald's operates its own delivery service called McDelivery in China and provides a 30-minute delivery guarantee to customers, which accounts for one-fifth of McDonald's revenue in mainland China \citep{mc2019, mc2020}. Grocery retailer chains such as Whole Foods Market, Co-op Food, Hema of Alibaba, 7Fresh Market of JD.com provide customers on-demand delivery services from their own stores within a few hours  \citep{hema2018, amzfresh2019}. Most recently, the COVID-19 pandemic has pushed more companies to expand their own delivery capacity and compete with platforms. According to a food delivery platform called Spread, one half of the restaurants on the platform make their own deliveries \citep{wsj2021b}.
Domino's Pizza, which is acclaimed for its delivery services, has delivered considerable sales growth with its own delivery capacity during the pandemic \citep{dmno2020}.

For on-demand delivery, achieving a high-quality delivery service in terms of both speed and reliability is critical. According to a 2019 national survey \citep{usfood2020}, cold food and delivery delays are the top two customer complaints for meal delivery services. Chain stores like Whole Foods and McDonald's are competing to deliver orders within the time-frame of hours or minutes. Satisfying such stringent on-time performance targets while maintaining a reasonable operational cost poses a big obstacle to the management of these delivery systems. 

Our study is motivated by a large grocery chain store in China. The store offers on-demand delivery services for grocery and prepared food (meal boxes). Each store serves a prespecified service region and delivers orders to customers using a dedicated fleet of drivers.  The store operates the system with multiple \textit{dispatch waves} (decision epochs): the planning horizon is divided into multiple time slots of equal lengths (15 minutes), so the orders placed in the same slot are bundled together and assigned to drivers who will be dispatched at a decision epoch. Each driver will be dispatched multiple times and perform multiple trips in the planning horizon. Given a limited fleet size, the company's goal is to optimize the overall on-time performance of delivery orders. 

Providing reliable on-time performance in last-mile delivery hinges on effective dispatching and routing of drivers. The studied delivery system features a highly dynamic and stochastic demand process, in which random customer orders arrive sequentially over a planning horizon. Customer locations and order quantities are both uncertain, and the operator (e.g., the store) can hardly preload orders or prespecify routes for the drivers in practice. As such, the operator needs to dispatch and route drivers dynamically in response to the realized customer orders. Specifically, due to a limited capacity, the system operator has to trade off the on-time performance of realized orders versus future orders. 

\subsection{Our Contributions}
How to dispatch and route a fleet of vehicles to fulfill random on-demand delivery orders quickly? Motivated by a large grocery chain store, we address this question by presenting a finite-horizon stochastic dynamic program for on-time delivery operations management. Our model captures the spatiotemporal heterogeneity and uncertainty of on-demand orders. Notably, because delivery drivers have to perform multiple trips within the planning horizon, we consider the interactions between dispatching and routing decisions explicitly.

Our key methodological contribution is a structured approximation framework that yields high-quality dispatching and routing decisions efficiently. Specifically, our framework estimates the cost-to-go function with a decomposed dispatching and routing policy. The estimation is then embedded into the dynamic program that outputs solutions in a rollout fashion. To this end, we integrate offline estimation and online rollout effectively. Our framework extends the existing approximate dynamic programming approaches in the vehicle routing literature to the multi-vehicle routing problem across multiple periods in a stochastic and dynamic setting. More importantly, we analyze the structural properties of our approximation framework and derive an approximation bound under large-demand scenarios. 

On the algorithmic side, we leverage the structure of the approximation model to develop computationally efficient algorithms by combining Benders decomposition and column generation, which allows an exact search of rollout policies. While a direct implementation with CPLEX fails to deliver solutions within an hour, the proposed decomposition algorithm finds optimal solutions in minutes. As a side product, our algorithm also leads to substantial improvement in solution times for an important class of vehicle routing problems against relevant state-of-the-art benchmarks.

We demonstrate the performance of our method on a real-world data set from our industry partner. Compared to the current solution policy of the company, our method yields 16\%-50\% improvement in delivery time. The improvement is further validated on a set of synthetic instances. From the empirical study, we quantify the value of dynamic dispatching and routing with different fleet sizes. We find that dynamic routing is more beneficial when the fleet size is not so large. We also discuss the value of increasing dispatch frequency, performing flexible order postponement, and varying the sample size under our framework, which leads to multiple prescriptions for further improving the on-time performance.

\subsection{Literature Review}\label{sec:review}
Our paper contributes to two streams of related literature: on-demand delivery operations and the vehicle routing literature focusing on dynamic problems.

\textbf{On-Demand Delivery Operations.}
Recently, on-demand delivery, particularly grocery and meal delivery, has received growing attention from transportation and operations management researchers. \cite{yildiz2019provably} solve the meal-delivery routing problem exactly with a simultaneous column- and row- generation, assuming perfect future information. They have performed extensive numerical experiments based on real-world data from Grubhub to validate the efficacy of their solutions. They highlight the importance of order bundling, driver shift scheduling, and demand management from the numerical study. In a relevant paper, \cite{reyes2018meal} propose optimization based algorithms and heuristics to solve the real-time assignment/dispatching problem in meal delivery. In contrast to our model, they do not capture the future order information in the assignment decisions. Nevertheless, we adopt similar metrics to measure the on-time performance of the delivery service. Based on a stylized queueing model, \cite{chen2020courier} analyze the optimal structure of the dispatching policy considering customers' patience level. They show that delivering multiple orders per trip is beneficial when the service area is large. In a general meal delivery context, \cite{ulmer2020restaurant} propose heuristic order assignment policies by introducing a time buffer cost as well as a postponement strategy. While they assume simplified assignment heuristics (not fully forward looking), our work aims to find assignment decisions that account for future order arrival uncertainties explicitly. \cite{liu2020time} study a meal delivery problem for a centralized kitchen and propose several ways to account for drivers' routing behaviors by integrating machine learning and optimization. They mainly focus on the single-period model, and only provide simple heuristics for the multiperiod setting.  Other aspects of on-demand delivery problems have also been studied, including the workforce scheduling \citep{ulmerworkforce}, supply management \citep{leidynamic} and demand management \citep{yildiz2020pricing}, and platform operations \citep{bahrami2021three}.


\textbf{Vehicle routing.}
The vehicle routing problem (VRP) has been a focal research topic of transportation and logistics since it was first proposed by \citet{dantzig1959truck}. 
According to the availability of information, the VRP can be classified into three basic variants, namely the static VRP, the stochastic VRP and the dynamic VRP. The static VRP has all input information available and all parameters in the problem are known and fixed. The stochastic VRP extends the static VRP by incorporating uncertain model parameters, including demands \citep{bertsimas1992vehicle}, travel time \citep{laporte1992vehicle, adulyasak2016models}, service times \citep{lei2012generalized}. The dynamic VRP, similar to the stochastic VRP, also has partial known input information when the routing plan is made, but the information is gradually revealed during the plan execution. The dynamism in most of the dynamic VRP originates from the online arrival of customer requests during the plan execution \citep{pillac2013review}. The driver dispatching and routing problem studied in this paper is a multiperiod problem, deciding routing plan to fulfill orders in the current period with an eye on the uncertain future orders. In terms of the single-period version of our problem, the most relevant static VRP is the multiple traveling repairman problem (MTRP) \citep{luo2014branch} whose objective to minimize the total arrival time at the customers. The MTRP has been tackled by various solution approaches, including mixed integer programming (MIP) \citep{nucamendi2016mixed,onder2017new}, branch-and-price \citep{luo2014branch}, and branch-and-cut \citep{muritiba2021branch}.

Among the dynamic VRP literature, the papers that are closest to our setting are \cite{azi2010exact, azi2012dynamic}, where the authors study the VRP with multiple delivery routes in a deterministic and stochastic context, respectively. In the stochastic setting, \cite{azi2012dynamic} develop a simulation based sample-scenario method combined with insertion and neighborhood search heuristics. In their paper, the main goal is to maximize the expected profits with the order acceptance decision, which is suitable for the same-day delivery environment. Our paper is focused on improving the on-time performance as highlighted by the emerging meal and grocery delivery services, where individual order rejection is not encouraged. In terms of methodologies, our paper is based on lookahead approximations in which the cost-to-go function is approximated by simple dispatching and routing policies (also called rollout policies, see \citealt{powell2019unified} for a detailed introduction). In contrast to existing lookahead methods that rely on heuristics to search for rollout policies in a restricted decision space \citep{cortes2009hybrid, goodson2013rollout, goodson2016restocking}, our approach integrates the rollout policy search and the decision making for the current state in one mixed integer linear program (MILP) and exploits its structure to enable exact rollout policy search efficiently.

When dispatching decisions are made at fixed intervals, \cite{klapp2018dynamic, klapp2018one} study a dynamic dispatch waves problem where a single vehicle is dispatched to serve orders on a network and on a one-dimensional line, respectively. They propose the a priori policy and several dynamic heuristic policies to solve the problem and show that dynamic policies can boost the system performance significantly. As their results only hold for the single-vehicle case, we demonstrate in this paper a framework to handle the general mutli-vehicle dispatching and routing problems with demand uncertainty. Our framework preserves preferable structural properties of the original problem, yielding a worst-case performance guarantee. On a high level, our proposed algorithms operationalize the batching policy proposed in \cite{bertsimas1993stochastic}.


\cite{voccia2019same} propose the same-day delivery problem (SDDP) that shares a similar structure to ours. While their objective is to maximize the expected number of fulfilled orders, we aim to minimize the expected delivery time, as motivated by our application in on-time delivery. 
 Because a complicated team orienteering problem has to be solved for every possible scenario, \cite{voccia2019same} apply neighborhood search heuristics as a solution subroutine, of which the optimality can be hardly guaranteed. In contrast, our decomposition-based algorithm integrates offline estimation and online rollout in a tractable manner. Note that our use of offline estimation is different from the offline-online approximate dynamic programming approach (ADP) proposed by \cite{ulmer2019offline}. Specifically, we do not require policy iterations to estimate and evaluate approximate policies for value function approximations. Notably, we extend their work on single-vehicle dynamic routing to the multi-vehicle setting with random demand, where a driver can take multiple trips in the planning horizon, and the need for coordination between vehicles across periods is prominent. Extensions to allow preemptive returns of vehicles and dynamic pricing of delivery deadline options are explored in \cite{ulmer2019preemptive} and \cite{ulmer2020dynamic}, respectively. We do not consider preemptive returns due to its implementation difficulties in the on-time delivery setting. 
 We refer interested readers to \cite{ulmer2019modeling} for an excellent review of relevant dynamic VRP papers. \cite{ulmer2019modeling} advocate the use of route-based models to bridge the gap between real-world applications and solution methodologies. Following a similar paradigm, our model has designed the route plan for realized orders in each epoch and specified the dispatching plan for future orders. 

\section{Problem Background and Description} \label{sec:background}
The studied on-time delivery problem is motivated by a large grocery chain store in China. The grocery chain operates in multiple cities across the country and adopts an omnichannel business model. In addition to serving in-store customers, each store offers on-demand delivery services to customers who place orders in a prespecified service region centered around the store. The delivery services cover a variety of products, from grocery goods to prepared meal boxes. Due to the high volume of demand, the company operates a separate channel for meal box delivery. Targeting stringent and reliable on-time performance, the company has hired a dedicated fleet of drivers to fulfill on-demand delivery tasks.\footnote{The use of dedicated drivers (in-house drivers) is not uncommon even for delivery platforms. Based on our communications with a leading meal delivery platform, a fleet of dedicated drivers can be deployed to serve high-demand restaurants.}

\begin{figure}[htbp]
	\subfloat[Spatial distribution\label{sfig:uniformt2}]{
		\includegraphics[width=0.40\linewidth]{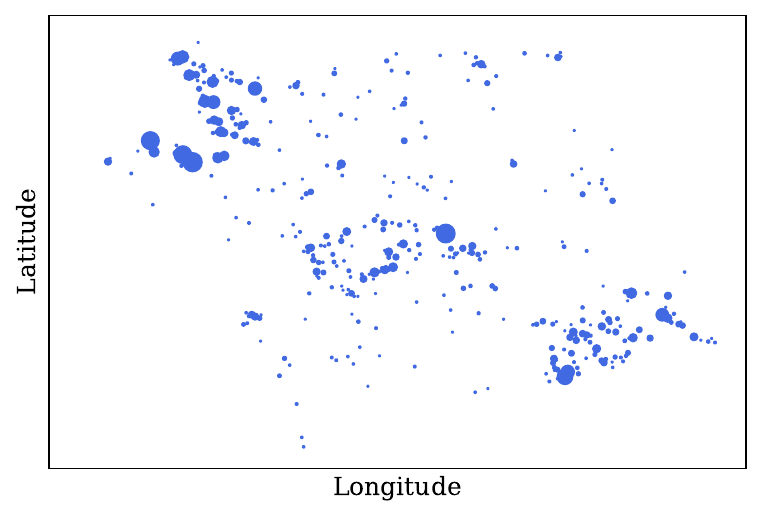}	}\hfill
	\subfloat[Temporal distribution\label{sfig:uniformwt2}]{
		\includegraphics[width=0.42\linewidth]{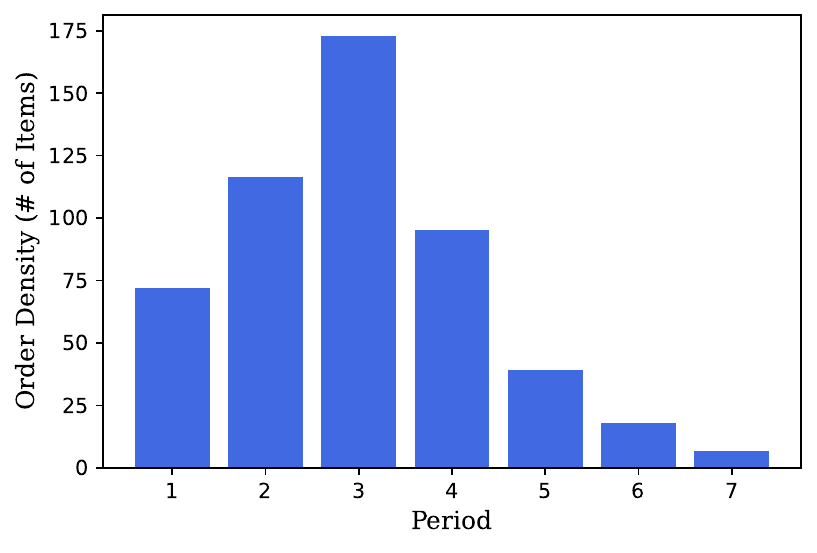}
	}\hfill
	\caption{Order Density Distribution}\label{fig:order}
\end{figure}

The delivery system operator of the company has specified a sequence of cutoff times to bundle customer orders together, corresponding to a set of dispatch waves. During the lunch peak hours, the cutoff times are [10 am, 10:15 am, 10:30 am, $\dots$, 11:45 am], making up seven 15-minute time slots (periods). As shown in Figure \ref{fig:order}, order density is  spatially and temporally heterogeneous, and there is a single demand peak in period 4. The operator processes orders in a batch process: the orders placed in the same time slot form a \textit{batch}, sharing the same delivery time target. For instance, the orders placed between 10:00 am and 10:15 am are promised to be delivered by 11:30 am. Once a batch of orders is collected, the store starts preparing the orders, and the operator will assign the batch of orders to available drivers and specify their routes (a visualization of this process is presented in Appendix \ref{sec:visual}). After the orders have been prepared (order preparation takes around 20 minutes), dispatched drivers will pick up orders at the store and perform deliveries. Typically, a driver can deliver multiple orders per trip (in many cases, more than five), which can take between 20 and 50 minutes.\footnote{Delivery boxes installed on the vehicles can maintain the freshness of orders during delivery.} Drivers will return to the store after finishing the assigned deliveries and become available for next dispatch.\footnote{The company does not allow preemptive returns of drivers for two reasons: (a) the online app allows customers to track the delivery process in detail, so having preemptive returns may cause customer confusion and complaints; (b) making preemptive returns may also give rise to fairness and equity concerns from customers.}  Because customers highly value delivery speed and promise reliability, the company desires a good dispatching and routing policy to minimize the delivery time while controlling their delivery fleet size and labor cost.\footnote{Because driver wage is paid based on the work duration (base payment) and the number of delivered orders (bonus payment), we do not consider driver travel cost directly. Nevertheless, our model can integrate travel cost into the objective function.}

\section{Driver Dispatching and Routing Model for On-Time Delivery}
In this section, we present the finite-horizon driver dispatching and routing model with multiple dispatch waves. A table that summarizes the notations used throughout the paper is included in Appendix \ref{sec:notations}. The delivery service system is operated for $N$ periods, with $\Delta$ denoting the length of each period, i.e., the time interval between two consecutive decision epochs. Our modeling framework does not require $\Delta$ to be stationary, but we assume the decision epochs are prespecified.  
We focus on the intraday operations with known shift schedule of drivers, i.e., the number of available drivers in period $n\in\{1,\dots, N\}$,  $\bar{K}^n$, is known and fixed. Each driver has a constant travel speed $v$ and can deliver at most $Q$ items per trip. For ease of discussion, we assume $\bar{K}^n = \bar{K}$ for $n=1,\dots, N$.

Denote the set of potential customer locations by $\mathcal{I}=\{1,\dots, I\}$.  At the beginning of time period $n$ (decision epoch $t_n$), we observe the number of orders realized between $[t_{n-1}, t_n)$, where $t_0$ denotes the start time of the service. These orders are to be assigned at $t_n$ and are characterized by their locations $\mathcal{I}^n\subseteq \mathcal{I}$ and quantities (number of ordered items) $\BFq^n = (q^n_1,\dots, q^n_I)\in \mathbb{N}^I$ ($q^n_{i} = 0$ if $i\notin \mathcal{I}^n$). The order quantities are integral with finite support ($q^n_{i}\leq Q$ without loss of generality). We assume a constant time $t_p$ for preparing and packaging the orders, and all the orders placed in period $n-1$ will be ready for delivery on $t_n + t_p$. Because the parameter $t_p$ can be estimated from data,  it is assumed to be known to the operator. We hereafter assume $t_p=0$ for ease of exposition,  and the incorporation of a positive $t_p$ is straightforward.

At decision epoch $t_n$, the operator has perfect information about which drivers are available for dispatch so they can pick up the orders at $t_n$.  
This is often the case for today's delivery system because drivers' smart phones are sending their real-time location information to the operator. We denote the driver status vector by $\BFzeta^n=(\zeta^n_1, \dots, \zeta^n_N)$. Specifically, $\zeta^n_{n'}\in \mathbb{N}$ is the number of en route drivers in period $n'$ due to the dispatching decisions made prior to period $n$. Note that only the $\zeta_{n'}^n$'s with $n'\geq n$ are meaningful, and we maintain the whole $\BFzeta^n$ for ease of reference. Summing the above information up, the state of the system at epoch $t_n$ is represented by $(\mathcal I^n, \BFq^n,\BFzeta^n)$.

The decision at the beginning of period $n$ is twofold: (1) we need to decide how many drivers to dispatch for the realized orders, which is denoted by variable $K^n$; (2) in the meanwhile, we also make the assignment of realized orders to the $K^n$ drivers and plan their routes. The dispatching decision echoes the scheduling decision, while the routing decision provides detailed execution plans. The dispatched drivers will return to the depot and become available again after finishing the assigned tasks. Following the company's practice, we assume all the orders in $\mathcal{I}^n$ are assigned to available drivers at $t_n$, i.e., orders placed between $[t_{n-1}, t_n)$ will not be assigned later than $t_n$. Such practice is preferable to reduce the wait time of orders at the store. Although allowing flexible order postponement can be beneficial, the additional gain may not be significant when the driver shift schedule and dispatching decisions are well optimized, as we numerically demonstrate in Section \ref{subsec:postpone}.  

We proceed to present the dynamic programming formulation. Let location $0$ be the depot where drivers are initially deployed and $\mathbb{Q}^{n}$ be the joint distribution of the customer locations and order quantities in period $n$ (for oders placed between $t_{n-1}$ and $t_n$). The system operator makes the joint dispatching and routing decision $\BFY^n = \{y_{ijk}^n\in\{0,1\}: i,j\in\mathcal{I}^n \cup \{0\}, k=1,\dots, \bar{K}\}$, where $y_{ijk}^n=1$ if driver $k$ is routed from $i$ to $j$ in period $n$ and 0 otherwise (note that the trip from $i$ to $j$ is not necessarily completed in period $n$). $y_{00k}^n=1$ indicates driver $k$ is not dispatched and stays at the depot.  The on-time performance measure for customer $i$ is denoted by $u_i(\BFY^n)$, which indicates the duration from the time the order is ready for dispatch until it is delivered following decision $\BFY^n$. Additionally, there is a hard delivery time target $L_{\max}$ for every order. The set of feasible decisions $\mathscr{D}(\mathcal I^n, \BFq^n,\BFzeta^n)$ must satisfy
\begin{align}
&y_{ijk}^n = 0, &\quad \forall i,j\in\mathcal{I}^n, k=1,\dots, \zeta^n_n, \label{const:avail} \\
&\sum_{i \in \mathcal{I}^n \cup \{0\}}y^n_{0ik} = \sum_{i \in \mathcal{I}^n \cup \{0\}}y^n_{i0k} = 1, &\quad \forall k = \zeta^n_n + 1,\ldots, \bar{K}, \label{const:depot}\\
&\sum_{j \in \mathcal{I}^n \cup \{0\}}y^n_{ijk} = \sum_{j \in \mathcal{I}^n \cup \{0\}}y^n_{jik}, &\quad \forall i \in\mathcal{I}^n, k = \zeta^n_n + 1,\ldots, \bar{K}, \label{const:flow}\\
&\sum_{i \in \mathcal{S}}\sum_{j \in \mathcal{S}}y^n_{ijk} \leq |\mathcal{S}| - 1, &\quad \forall \mathcal{S} \subseteq \mathcal{I}^n, k = \zeta^n_n + 1,\ldots, \bar{K}, \label{const:sec} \\
&\sum_{i \in \mathcal{I}^n}\sum_{j \in \mathcal{I}^n \cup \{0\}}q^n_iy^n_{ijk} \leq Q, &\quad \forall k = \zeta^n_n + 1,\ldots, \bar{K}, \label{const:capacity}\\
& u_i(\BFY^n) \leq L_{\max}, &\quad \forall i\in\mathcal{I}^n,   \label{const:duration}
\end{align}
where constraints (\ref{const:avail}) impose the driver availability condition, i.e., the drivers who are occupied in period $n$ due to the assigned delivery tasks can not be dispatched in period $n$. Constraints (\ref{const:depot}) ensure that each driver trip must start from and end at the depot. Constraints (\ref{const:flow}) and (\ref{const:sec}) are the flow conservation constraints and the subtour elimination constraints, respectively. Constraints (\ref{const:capacity}) ensure driver capacity is not violated and constraints (\ref{const:duration}) respect the hard delivery time target (for brevity we move the detailed representation of $u_i(\BFY^n)$ to Appendix \ref{sec:hdd}).

The objective is to minimize the total expected delivery time of orders in the planning horizon.  Let $l_k^n(\BFY^n)$ denote the route duration (including both travel time and service time) of driver $k$ dispatched in period $n$. The finite-horizon stochastic dynamic program for on-time delivery can be formulated with the value (cost-to-go) functions $\mathcal{H}_n(\mathcal I^n, \BFq^n,\BFzeta^n)$ as
\begin{align}
&\mathcal{H}_n(\mathcal I^n, \BFq^n,\BFzeta^n) = \min_{\BFY^n\in \mathscr{D}(\mathcal I^n, \BFq^n,\BFzeta^n)}\left\{    \sum_{i \in \mathcal I^n} u_i(\BFY^n) + \mathbb{E}_{\mathbb{Q}^{n+1}} \left[\mathcal{H}_{n+1}(\mathcal I^{n+1},  \BFq^{n+1}, \BFzeta^{n+1}) \right] \right\}, \label{eqn:bell1} \\
&\mathcal{H}_N(\mathcal I^N, \BFq^N, \BFzeta^N) = \min_{\ \BFY^N \in \mathscr{D}(\mathcal I^N, \BFq^N, \BFzeta^N)}\left\{  \sum_{i \in \mathcal I^N} u_i(\BFY^N)\right\},
&
\end{align}
with the transition constraints for driver availability:
\begin{align}
&\zeta_{n'}^{n+1}= \zeta_{n'}^n + \sum_{k=1}^{\bar{K}} \mathds{1}(l^n_k(\BFY^n) > t_{n'} - t_n),\quad \forall n'=n+1, \dots, N,\ n = 1, \dots, N \label{eqn:transition} ,
\end{align}
where $\mathds{1}(l^n_k(\BFY^n) > t_{n'} - t_n)$ is an indicator variable that equals 1 if driver $k$ can not return to the depot before period $n'$ given decision $\BFY^n$. Note that the choice of the on-time performance measure is flexible, and our model can incorporate other metrics such as ready-to-door time and click-to-door time overage.  We refer to the above dynamic program as \textbf{JDR}.

Due to the capacity and delivery time constraints, the dynamic program may not always be feasible when the number of available drivers ($\bar{K}$) is small. As we will discuss later, even when there is an adequate driver schedule, a smart dispatching policy is necessary to yield a feasible solution for every period. In practice, we can introduce simple recourse rules to tackle infeasible scenarios, such as calling additional drivers from third-party platforms. We will discuss this option in Section \ref{sec:sim}.


\section{A Structured Approximation Approach}
Because both the state space and the action space are high dimensional, 
\textbf{JDR} can not be solved exactly.  Even when the demand is deterministic, the resulting multiperiod dispatching and routing problem is NP-hard and potentially time consuming to solve \citep{klapp2018dynamic}. The combinatorial nature of the problem and the complicated dependence on the random demand stresses the difficulty of analysis and optimization.
Therefore, it is not uncommon to see companies use simple myopic policies to dispatch and route drivers in delivery planning: the dispatching and routing decisions are obtained to optimize the on-time performance of the current batch of orders without accounting for future order arrivals. However, in the considered planning horizon, a driver must perform multiple trips and, thus, travel back and forth between the store and customers (all the orders must be first picked up at the store). The dispatching and routing decision made for the current batch will decide the driver availability in the future periods (as shown in Equation \eqref{eqn:transition}). Ignoring this interaction can severely exacerbate the long-run system performance, e.g., when drivers are sent out blindly to serve realized orders, and none of them are available for the next dispatch wave. A forward-looking dispatching and routing policy is desired to properly trade off the delivery time of realized orders versus future orders.



To yield high-quality solutions in real time, we develop a tractable approximation framework for the studied stochastic dynamic program. At a high level, our framework estimates the cost-to-go function through a parameterized dispatching and routing policy that combines myopic routing with anticipatory dispatching. 
The estimated cost-to-go function will then help identify the best dispatching and routing decision for the current state. In contrast to existing value approximation methods, we show that our approximation framework preserves structural properties of the true cost-to-go function, which helps bound the approximation ratio.

The key to establishing the approximate cost-to-go function is modeling the impact of the decision (or post-decision state) on future costs.  The dispatching and routing decision affects the future delivery costs through restricting the number of available drivers in the remaining planning horizon. 
Specifically, when more drivers are dispatched for the current period, fewer drivers will be available for delivery in the following periods. Similarly, when drivers are assigned longer routes, future delivery capacity will be affected because it takes a longer time for the dispatched drivers to return to the depot. The timing of dispatch waves should be respected so that the drivers' availability information can be accounted for properly. To capture this delicate relationship between future driver supply and delivery cost, we approximate the cost-to-go function by the sum of single-period value functions under myopic routing policies. Specifically, let  $\mathcal H^s(K^n,\mathcal{I}^n,\mathbf q^n)$ denote the single-period optimal delivery cost with $K^n$ dispatched drivers when the realized customer locations and order quantities are $\mathcal{I}^n$ and $\BFq^{n}$, respectively. Denote $\omega_{m}^{n'}(K^m)$ by the number of en route drivers in period $n'$ out of the $K^m$ drivers dispatched in period $m$. Then the expected cost-to-go function $\mathbb{E}_{\mathbb{Q}^{n}} \left[\mathcal{H}_{n}(\mathcal I^{n},  \BFq^{n}, \BFzeta^{n}) \right]$ is approximated by
\begin{align}
    \textbf{APT$^n(\BFzeta^{n})$:} \quad & \min_{K^{n'} \in \mathbb{N}} \sum_{n'=n}^N   \mathbb{E}_{\mathbb Q^{n'}} \left[\mathcal H^s(K^{n'},\mathcal{I}^{n'},\BFq^{n'})\right] \nonumber \\
    & s.t.\ \sum_{m=n}^{n'} \omega_{m}^{n'}(K^m)  \leq \bar{K} - \zeta^n_{n'},\quad \forall n'=n,\dots, N, \label{const:availability}
\end{align}
where $\mathbb{E}_{\mathbb Q^{n}}\left[\mathcal H^s(K^n,\mathcal{I}^{n},\BFq^{n}) \right]$ is the expected single-period optimal delivery cost, summing over all possible realizations of $\mathcal{I}^{n}$ and $\BFq^{n}$.We can estimate it by offline simulations based on historical data or a fitted probability distribution: $ \mathbb{E}_{\mathbb Q^{n}}\left[\mathcal H^s(K^n,\mathcal{I}^{n},\BFq^{n}) \right] = \sum_{h=1}^H \mathcal H^s(K^n,\mathcal{I}^n_h,\mathbf q^n_h)/H$ for $H$ samples of customer locations and orders (although the possible scenarios can be many, a finite sample of historical data can capture the general spatiotemporal pattern of demand). Constraints (\ref{const:availability}) ensure the number of dispatched and en route drivers does not exceed $\bar{K}$ in every period.
Note that a dispatched driver's en route time is at least one period (i.e., a driver can not be dispatched again within a period), so $\omega^{n'}_{n'}(K^{n'}) = K^{n'}$ for $n'=n,\dots, N$. However, for $n'>m$, $\omega_{m}^{n'}(K^m)$ is uncertain due to the stochastic nature of demand, and we treat it as a parameter that can be calibrated or tuned from offline simulations.

\begin{figure}[htbp]
	\centering
	\includegraphics[width=0.6\linewidth]{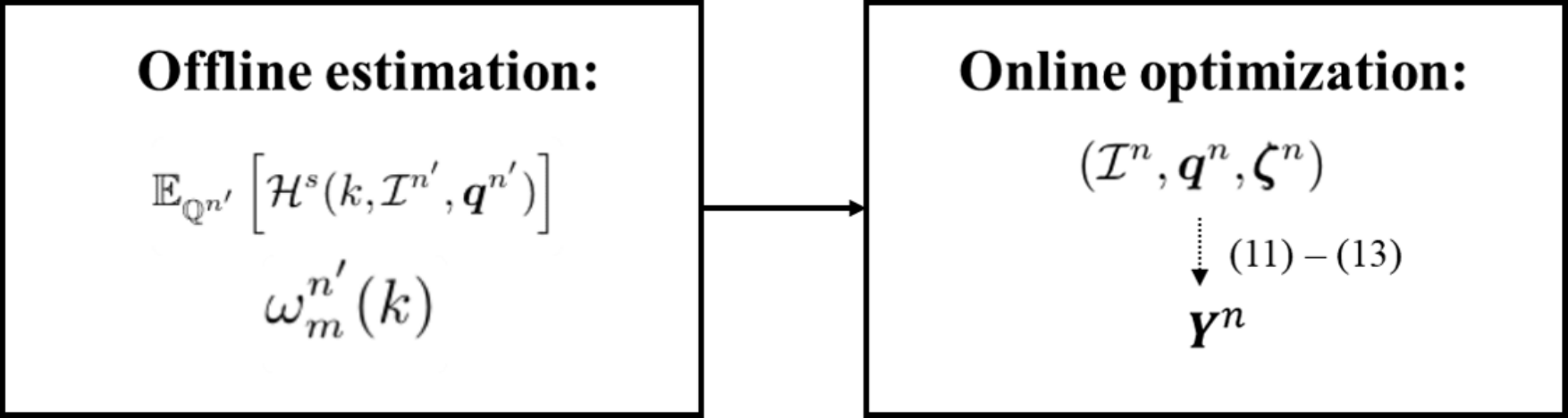}\hfill
	\caption{An Illustration of AJRP}\label{fig:illustration}
\end{figure}

The above approximation scheme estimates the expected cost-to-go function by a decomposed dispatching and (myopic) routing heuristic. It can be viewed as a stochastic lookahead approach based on rollout policies in approximate dynamic programming (the readers may find a detailed introduction to lookahead methods in \citealt{powell2011approximate}). Under this lookahead approach, the routing of future orders is assumed to be myopic when evaluating $\mathcal H^s(K^n,\mathcal{I}^n,\mathbf q^n)$. Albeit myopic in routing for each period, it strives to capture the relationship between driver supply and delivery cost through detailed modeling of dispatching with respect to dispatch waves. 
Note that the heuristic myopic policies (rollout policies) will not be implemented but only to facilitate the decision selection in the current decision epoch (so we do not need to foresee all possible future scenarios). 
Specifically, the approximation \textbf{APT}$^{n+1}(\BFzeta^{n+1})$ is used in dynamic program (\ref{eqn:bell1}) to find the dispatching and routing decision at decision epoch $n$ and state $(\mathcal I^n, \BFq^n,\BFzeta^n)$:
\begin{align*}
\min_{\BFY^n\in \mathscr{D}(\mathcal I^n, \BFq^n,\BFzeta^n)}\ & \left\{    \sum_{i \in \mathcal I^n} u_i(\BFY^n) +\min_{K^{n'} \in \mathbb{N}} \sum_{n'=n+1}^N   \mathbb{E}_{\mathbb Q^{n'}} \left[\mathcal H^s(K^{n'},\mathcal{I}^{n'},\BFq^{n'})\right] \right\} \\
s.t.\ &  \sum_{m=n+1}^{n'} \omega_{m}^{n'}(K^m)  \leq \bar{K} - \zeta^{n+1}_{n'},\quad \forall n'=n+1,\dots, N.
\end{align*}
Introducing binary variables $x_k^{n'}$ to indicate if $k$ drivers are dispatched in period $n'$ ($x_k^{n'}=1$), and leveraging the state transition equation (\ref{eqn:transition}), the above program can be rewritten as
\begin{align}
\min_{\substack{\BFY^n\in \mathscr{D}(\mathcal I^n, \BFq^n,\BFzeta^n)\\ x_k^{n'}\in\{0,1\}}} \ &     \sum_{i \in \mathcal I^n} u_i(\BFY^n) + \sum_{n'=n + 1}^N \sum_{k=0}^{\bar{K}} x_k^{n'} \mathbb{E}_{\mathbb Q^{n'}}\left[\mathcal H^s(k,\mathcal{I}^{n'},\BFq^{n'})  \right], \label{ajrp:obj}\\
s.t. \ & \sum_{m=n+1}^{n'} \sum_{k=0}^{\bar{K}} \omega_{m}^{n'}(k) x_k^{m}   \leq \bar{K} - \zeta^n_{n'} - \sum_{k=1}^{\bar{K}} \mathds{1}(l^n_k(\BFY^n) > t_{n'} - t_n) ,\quad \forall n'=n+1,\dots, N, \label{ajrp:const1}\\
&\sum_{k=0}^{\bar{K}} x_k^{n'} =1,\quad \forall n' = n+1,\dots, N, \label{ajrp:const2}
\end{align}
where constraints (\ref{ajrp:const2}) ensure the number of dispatched drivers in every period can only take an integral value between 0 and $\bar{K}$. We refer to the resulting approximate joint dispatching and routing policy as AJRP. We illustrate how AJRP is solved by combining offline estimation and online rollout in Figure \ref{fig:illustration}. The rollout policy is parameterized by the single-period cost functions $ \{\mathbb{E}_{\mathbb Q^{n}}\left[\mathcal H^s(k,\mathcal{I}^{n},\BFq^{n}) \right]\}_{\forall (n,k)}$ and the driver state functions $\{\omega_{m}^{n'}(k)\}_{\forall (m,n',k)}$. In order to enumerate all possible dispatching decisions, we evaluate the single-period cost functions for all feasible integer values of $k$ in $[0, \bar{K}]$.

We now describe structural properties of our approximate cost-to-go function and provide a bound on the approximation ratio. First, we show that the optimal objective value of  \textbf{APT}$^n(\BFzeta^{n})$, denoted by $V_{APT}^n(\BFzeta^{n})$, is increasing in $\BFzeta^{n}$ for any choices of $\omega_{m}^{n'}(K^m)\geq 0$, which is consistent with $\mathcal{H}_n(\mathcal I^n, \BFq^n,\BFzeta^n) $.
\begin{lemma} \label{lem:2}
	(i) $\mathcal{H}_n(\mathcal I^n, \BFq^n,\BFzeta_1^n) \geq \mathcal{H}_n(\mathcal I^n, \BFq^n,\BFzeta_2^n)$ for $\BFzeta_1^n \geq \BFzeta_2^n$; (ii) Given a set of nonnegative $\{\omega_{m}^{n'}(K^m)\}_{\forall (m,n')}$, $V_{APT}^n(\BFzeta_1^{n}) \geq V_{APT}^n(\BFzeta_2^{n})$ for $\BFzeta_1^n \geq \BFzeta_2^n$.
\end{lemma}
Therefore, our approximation scheme maintains the monotonicity property of the true value function. Next, as shown in the following theorem, the proposed approximation is exact for the last period and can provide lower and upper bounds of the expected cost-to-go function with appropriate values of $\omega_{m}^{n'}(K^m)$. Before stating the theorem, we introduce the class of static myopic policies $\{\pi^{sm}\}$, wherein the number of dispatched drivers in each period is state independent, and the routing decision is myopic, i.e., we route drivers in a way that only minimizes the single-period cost. Let $\textbf{JDR}^n(\BFzeta^{n})$ denote the joint dispatching and routing problem starting in period $n$ with driver status $\BFzeta^{n}$, after compressing the demand information.

\begin{theorem} \label{them:approx}
Under the assumption that there exists a feasible static myopic policy to $\textbf{JDR}^n(\BFzeta^{n})$, the approximation \textbf{APT}$^n(\BFzeta^{n})$ can serve as lower and upper bounding problems of $\mathbb{E}_{\mathbb{Q}^{n}} \left[\mathcal{H}_{n}(\mathcal I^{n},  \BFq^{n}, \BFzeta^{n}) \right]$ with appropriate choices of $\{\omega_{m}^{n'}(K^m)\}_{\forall (m,n')}$. Furthermore, this approximation is exact for the last period.
\end{theorem}

The assumption of Theorem \ref{them:approx} will be satisfied when the fleet size is not too small relative to $\BFzeta^{n}$, otherwise any static myopic policy is infeasible, and a feasible policy must be fully adaptive to the realization of $(\mathcal I^{n},  \BFq^{n}, \BFzeta^{n})$. Nevertheless, the nonexistence of a feasible static myopic policy does not exclude the feasibility of problem \textbf{APT}$^n(\BFzeta^{n})$, which can still be solved to obtain a reasonable approximation to the value function.  Theorem \ref{them:approx} implies that AJRP is optimal for $N=2$.

\begin{corollary} \label{lem:n2}
	AJRP is optimal for \textbf{JDR} when $N=2$.
\end{corollary}

As indicated by Theorem \ref{them:approx}, the choice of $\{\omega_{m}^{n'}(K^m)\}_{\forall (m,n')}$ steers the relationship between \textbf{APT}$^n(\BFzeta^{n})$  and the true cost-to-go function.  Recall that $\omega_{m}^{n'}(K^m)$ reflects the number of en route drivers in period $n'$ out of the $K^m$ drivers who are dispatched in period $m$. Hence, we can evaluate $\omega_m^{n'}(K^m)$ by offline simulation using myopic routing policies. As such, the evaluation of $ \{\mathbb{E}_{\mathbb Q^{n}}\left[\mathcal H^s(k,\mathcal{I}^{n},\BFq^{n}) \right]\}_{\forall (n,k)}$ and $\{\omega_{m}^{n'}(K^m)\}_{\forall (m,n')}$ can be performed simultaneously. Let $\bar{\omega}_{m}^{n'}(K^m)$ denote the estimated average value of $\omega_m^{n'}(K^m)$ from simulation and $\bar{V}_{APT}^n(\BFzeta^{n})$ denote the optimal objective value of \textbf{APT}$^n(\BFzeta^{n})$ with the choice of $\{\bar{\omega}_{m}^{n'}(K^m)\}_{\forall (m,n')}$.
The following proposition establishes the relationship between $\bar{V}_{APT}^n(\BFzeta^{n})$ and $\mathbb{E}_{\mathbb{Q}^{n}} \left[\mathcal{H}_{n}(\mathcal I^{n},  \BFq^{n}, \BFzeta^{n}) \right]$.

\begin{proposition} \label{prop:approx}
Under the assumption that there exists a feasible static myopic policy to $\textbf{JDR}^n(\BFzeta^{n})$, $\bar{V}_{APT}^n(\BFzeta^{n})$ is finite, and there exists an instance specific $\vartheta>0$ such that
\[ 1/\vartheta \leq  \frac{\bar{V}_{APT}^n(\BFzeta^{n})}{\mathbb{E}_{\mathbb{Q}^{n}} \left[\mathcal{H}_{n}(\mathcal I^{n},  \BFq^{n}, \BFzeta^{n}) \right]} \leq  \vartheta.\]
Furthermore, there exists a positive constant $M$ such that  $\bar{V}_{APT}^n(\BFzeta^{n}) = \mathbb{E}_{\mathbb{Q}^{n}} \left[\mathcal{H}_{n}(\mathcal I^{n},  \BFq^{n}, \BFzeta^{n}) \right]$ when $\bar{K} \geq M$.
\end{proposition}

Proposition \ref{prop:approx} shows that the ratio of the approximation value $\bar{V}_{APT}^n(\BFzeta^{n})$ and the true value $\mathbb{E}_{\mathbb{Q}^{n}} \left[\mathcal{H}_{n}(\mathcal I^{n},  \BFq^{n}, \BFzeta^{n}) \right]$ can be bounded, which implies that AJRP has a bounded approximation ratio. 
As the driver pool becomes sufficiently large, the proposed approximation policy using $\bar{V}_{APT}^n(\BFzeta^{n})$ is optimal.
Although the approximation ratio is instance-dependent, we leverage the above structural results to prove a worst-case performance guarantee under large-demand scenarios. Without loss of generality, we assume demand locations are uniformly distributed in a bounded Euclidean service region of area $A$. Let $\bar{r}$ denote the average travel distance from the depot to a customer in the service region and $s$ denote the on-site service time of each order.

\begin{theorem} \label{prop:bound}
	Assuming there are at least $I^*$ realized customer locations in each period, and each customer orders exactly one item, the approximation ratio $\vartheta$ satisfies that for large $I^*$,
	\begin{align*}
	\vartheta \lessapprox \frac{\bar{r}/v + (Q+1)s/2 + \beta(Q-1)\sqrt{A}/(2v\sqrt{I^*})}{\bar{r}/v + s},
	\end{align*}
	where $\beta$ is a constant.
\end{theorem}

The above result bounds the approximation ratio of AJRP for systems with large demand, where we utilize the asymptotic analysis of the TSP tour length \citep{beardwood1959shortest, steele1981subadditive}. Based on \cite{applegate2010using},  the constant satisfies $0.6250 \leq \beta \leq 0.9204$. The derived bound depends on the geometry of the service region through $\bar{r}$ and $A$. Intuitively, the problem facing a smaller capacity $Q$ will result in a tighter bound because there is less room for dispatching and routing optimization. For a practical case where $\bar{r}/v = 15$ minutes, $s=2$ minutes, and $Q=10$, the computed upper bound is approximately 1.53 when $I^*$ is large. The assumption of a uniform demand distribution is not critical, and the analysis can be extended to general demand distribution functions. 

In the dispatching and routing literature, the commonly used heuristics and value function approximation methods do not enjoy performance guarantees. Theorem \ref{prop:bound} gives a characterization of the approximation ratio of AJRP under certain circumstances and bounds the performance gap. 
Moreover, our approximation enables a computationally efficient solution framework. The single-period cost functions $ \{\mathbb{E}_{\mathbb Q^{n}}\left[\mathcal H^s(k,\mathcal{I}^{n},\BFq^{n}) \right]\}_{\forall (n,k)}$ can be evaluated offline, which is facilitated by a specialized single-period optimization algorithm detailed in Section \ref{sec:rea}. In particular, the decomposable structure of AJRP gives rise to a Benders decomposition solution approach that admits verifiably optimal solutions quickly.

\section{A Logic Benders Decomposition Based Solution Framework}
\begin{figure}[h]
	\centering\scalebox{0.4}{\includegraphics{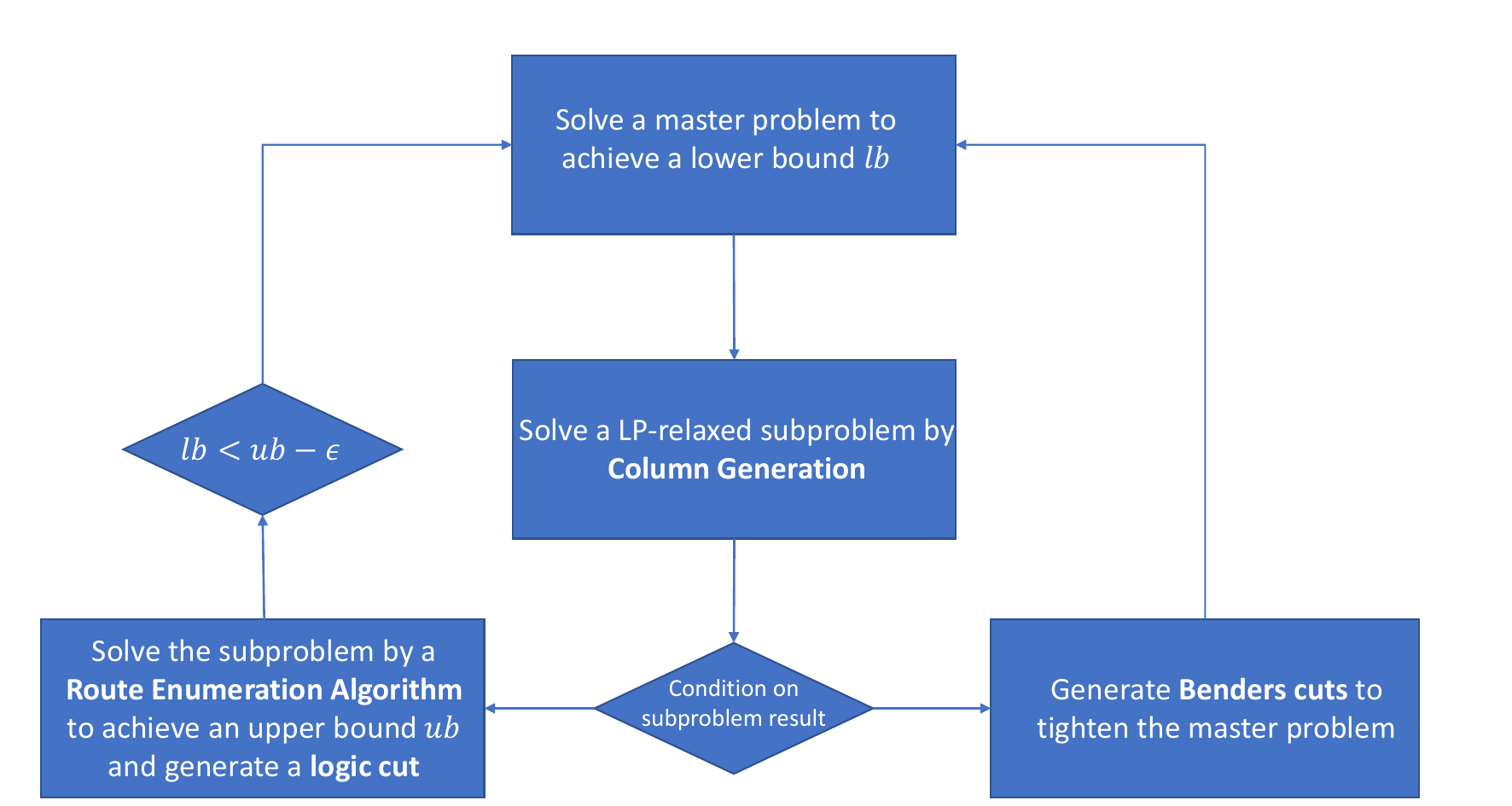}}\\
	\centering\caption{The Proposed Solution Framework Based on Benders Decomposition}
	\label{fig:benders}
\end{figure}

Although the AJRP formulation can be tackled by off-the-shelf solvers such as CPLEX and Gurobi, the solution time is often a bottleneck to practical real-time implementation. According to our preliminary computational experiments, a direct implementation of the AJRP formulation in CPLEX can not deliver optimal solutions in one hour, even for the smallest instances. In this section, we develop an efficient solution framework to obtain verifiably optimal solutions by exploiting the structure induced by AJRP.  Specifically, the driver dispatching and routing decisions under AJRP can be organized in a two-stage manner, i.e., the number of dispatched drivers in the first stage and the detailed routing plan for each driver in the second stage. Based on this observation, we propose a logic Benders decomposition method to solve AJRP. Figure \ref{fig:benders} provides an overview of our solution framework -- the proposed algorithm iteratively solves a master problem and a sub-problem until the optimality gap is small enough. In each iteration, the master problem is solved  to obtain a lower bound, and the LP relaxation of the sub-problem is solved by column generation. If the optimal LP cost of the sub-problem plus the cost of the master problem's solution is large enough to cut off the master problem's solution, a Benders cut is added to the master problem. Otherwise, the sub-problem is solved exactly to achieve an optimal integer solution and 
update the upper bound. Meanwhile, a logic Bender cut is added to the master problem to cut off the master problem's solution. In this section, we first introduce the Benders decomposition formulation and then describe the proposed column generation and route enumeration algorithms for solving the subproblems efficiently.

\subsection{Logic Benders Decomposition}
 We start by transforming arc based formulation (\ref{ajrp:obj}) -- (\ref{ajrp:const2}) to a route based formulation of AJRP, because the LP relaxation of the route based formulation provides a much better lower bound than the arc based formulation. Let $\mathcal{R}^n$ be the set of feasible routes for the orders in period $n = 1,\ldots,N$, $c_r$ be the cost of route $r \in \mathcal{R}^n$, $\alpha_{i,r}\in\{0,1\}$ indicate if location $i \in \mathcal{I}^n$ is served in route $r \in \mathcal{R}^n$, and $\beta_{n',r}\in\{0,1\}$ indicate if route $r \in \mathcal{R}^n$ is running in period $n' = n + 1, \ldots,N$. Further, we use binary variables $\theta_r$ to indicate whether route $r \in \mathcal{R}^n$ is assigned to a driver, and binary variables $z^{n'}_k$ to indicate whether there exist $k = 0,\ldots,\bar{K}$ drivers who are dispatched in period $n$ but are still occupied in period $n' = n + 1,\ldots,N$.  The resulting route-based formulation is
\begin{align}
\textbf{F1:} \quad \min &     \sum_{r \in \mathcal{R}^n} \theta_r c_r + \sum_{n'=n + 1}^N \sum_{k=0}^{\bar{K}} x_k^{n'} \mathbb{E}_{\mathbb Q^{n'}}\left[\mathcal H^s(k,\mathcal{I}^{n'},\BFq^{n'})  \right], \label{f1:obj}\\
s.t. \ & \sum_{m=n+1}^{n'} \sum_{k=0}^{\bar{K}} x_k^{m} \omega_{m}^{n'}(k)  \leq \bar{K} - \zeta^{n'}_{n}  - \sum_{k = 0}^{\bar{K}}kz^{n'}_k,\quad \forall n'=n+1,\dots, N, \label{f1:con1}\\
&\sum_{k=0}^{\bar{K}} x_k^{n'} =1,\quad \forall n' = n+1,\dots, N, \label{f1:con2}\\
&\sum_{k=0}^{\bar{K}} z_k^{n'} =1,\quad \forall n' = n+1,\dots, N, \label{f1:con9}\\
&\sum_{r \in \mathcal{R}^n}\beta_{n',r}\theta_r \leq \sum_{k = 0}^{\bar{K}}kz^{n'}_k,\quad \forall n'=n+1,\dots, N, \label{f1:con3}\\
&\sum_{r \in \mathcal{R}^n}\theta_r \leq \bar{K} - \zeta^n_n,   \label{f1:con4}\\
&\sum_{r \in \mathcal{R}^n} \alpha_{i,r}\theta_r =1, \quad \forall i\in \mathcal{I}^n, \label{f1:con5}\\
&x^{n'}_k, z^{n'}_k \in \{0,1\}, \quad \forall k = 0,\ldots,\bar{K},  n' = n + 1, \ldots,N, \label{f1:con6}\\
&\theta_r \in \{0,1\}, \quad \forall r \in \mathcal{R}^n. \label{f1:con8}
\end{align}
The objective function (\ref{f1:obj}) includes the cost of the current period and the approximate expected future cost. Constraints (\ref{f1:con1}) and (\ref{f1:con4}) ensure the number of occupied drivers in the future periods and the current period do not exceed the capacity (maximum number of available drivers), respectively. Constraints (\ref{f1:con2}) and (\ref{f1:con9}) enforce the convexity of variables $\mathbf{x}$ and $\mathbf{z}$, respectively. Constraints (\ref{f1:con3}) are the linking constraints between variables $\pmb{\theta}$ and $\mathbf{z}$. Constraints (\ref{f1:con5}) guarantee that each order of the current period is assigned to a driver.

Problem \textbf{F1} can be decomposed into a master problem that only involves dispatching decisions, and a subproblem consisting of the routing variables. Specifically, the master problem includes $(\mathbf{x}, \mathbf{z})$ and the subproblem decides on $\mathbf{\pmb{\theta}}$. The master problem is formulated as
\begin{align}
\textbf{MF1:} \quad \min &\sum_{n'=n + 1}^N \sum_{k=0}^{\bar{K}} x_k^{n'} \mathbb{E}_{\mathbb Q^{n'}}\left[\mathcal H^s(k,\mathcal{I}^{n'},\BFq^{n'})  \right], \label{f1:obj2}\\
s.t. \ & (\ref{f1:con1}), (\ref{f1:con2}), (\ref{f1:con9}), (\ref{f1:con6}). \nonumber
\end{align}
Given a feasible solution ($\bar{\textbf{x}}$,$\bar{\textbf{z}}$) of master problem \textbf{MF1}, the subproblem is a single-period routing problem as follows:
\begin{align}
\textbf{SF1:} \quad \min &  \sum_{r \in \mathcal{R}^n} \theta_r c_r, \label{sf1:obj}\\
s.t. \ &\sum_{r \in \mathcal{R}^n}\beta_{n',r}\theta_r \leq \sum_{k = 0}^{\bar{K}}k\bar{z}^{n'}_k,\quad \forall n'=n+1,\dots, N, \label{sf1:con1}\\
&(\ref{f1:con4}), (\ref{f1:con5}), (\ref{f1:con8}). \nonumber
\end{align}

Because the subproblem \textbf{SF1} is an integer program, we relax it as a linear program to derive the Benders cuts. The relaxed formulation \textbf{RF1} is
\begin{align}
\textbf{RF1:} \quad \min &\ (\ref{sf1:obj}) \nonumber\\
s.t. \ & (\ref{f1:con4}), (\ref{sf1:con1}) \nonumber\\
&\sum_{r \in \mathcal{R}^n} \alpha_{i,r}\theta_r \geq 1, \quad \forall i\in \mathcal I^n, \label{rf1:con1}\\
&\theta_r \geq 0, \quad \forall r \in \mathcal{R}^n. \label{rf1:con2}
\end{align}
Let $\mu_{n'}$ $(n' = n + 1,\ldots,N)$, $\mu_n$ and $\nu_i$ $(i \in \mathcal{I}^n)$ be the dual variables of constraints (\ref{sf1:con1}), (\ref{f1:con4}) and (\ref{rf1:con1}), respectively, and $\pmb{\Theta}$ and $\pmb{\Lambda}$ be the set of extreme points and extreme rays of problem \textbf{RF1}'s dual problem, respectively. We derive a relaxation of problem \textbf{F1} as
\begin{align}
\textbf{F2:} \quad \min &\sum_{n'=n + 1}^N \sum_{k=0}^{\bar{K}} x_k^{n'} \mathbb{E}_{\mathbb Q^{n'}}\left[\mathcal H^s(k,\mathcal{I}^{n'},\BFq^{n'})  \right] + \eta, \label{f2:obj}\\
s.t. \ &\sum_{n' = n + 1}^{N}\sum_{k = 0}^{\bar{K}}\mu_{n'}kz^{n'}_k + \mu_n(\bar{K} - \zeta^n_n) + \sum_{i \in \mathcal{I}^n}\nu_i \leq \eta, \quad \forall (\pmb{\mu},\pmb{\nu}) \in \pmb{\Theta} \label{f2:con1}\\
&\sum_{n' = n + 1}^{N}\sum_{k = 0}^{\bar{K}}\mu_{n'}kz^{n'}_k + \mu_n(\bar{K}  - \zeta^n_n) + \sum_{i \in \mathcal{I}^n}\nu_i \leq 0, \quad \forall (\pmb{\mu},\pmb{\nu}) \in \pmb{\Lambda} \label{f2:con2}\\
& (\ref{f1:con1}), (\ref{f1:con2}), (\ref{f1:con9}), (\ref{f1:con6}), \nonumber
\end{align}
where constraints (\ref{f2:con1}) and (\ref{f2:con2}) are the optimality Benders cuts and the infeasibility Benders cuts, respectively.

Note that the sizes of $\pmb{\Theta}$ and $\pmb{\Lambda}$ are exponential, so Benders cuts (\ref{f2:con1}) and (\ref{f2:con2}) cannot be enumerated beforehand. Instead, they are generated dynamically by solving problem \textbf{RF1}. Meanwhile, because problem \textbf{F2} involves only dispatching related decision variables, problem \textbf{SF1} has to be exactly solved to get the detailed routing plan. Therefore, the Benders decomposition solves the relaxed master problem \textbf{F2} and subproblems \textbf{SF1} and \textbf{RF1} successively. The implementation details of the Benders decomposition are presented in Appendix \ref{sec:dbd}.
Because the optimality cuts (\ref{f2:con1}) and the infeasibility cuts (\ref{f2:con2}) are derived from the LP relaxation of subproblem \textbf{SF1}, problem \textbf{F2} is a relaxation of problem \textbf{F1}. As a result, the dispatching decision obtained from the solution of problem \textbf{F2} may be infeasible or non-optimal for problem \textbf{F1}. Specifically, if subproblem \textbf{SF1} is infeasible, then the optimal solution of problem \textbf{F2} is also infeasible for problem \textbf{F1}. If the cost of problem \textbf{F2}'s optimal solution plus the cost of subproblem \textbf{SF1}'s solution is larger than the cost of the current best solution of problem \textbf{F1}, then the optimal solution of problem \textbf{F2} is non-optimal with respect to problem \textbf{F1}. Suppose $(\pmb{\bar{x}}, \pmb{\bar{z}})$  is such a solution, the following logic Benders cut is added to problem \textbf{F2} to cut it off:
\begin{align}
  &\sum_{n' = n + 1}^{N}\sum_{k = 0}^{\bar{K}}\left\{\mathds{1}(\bar{x}^{n'}_k = 0)x^{n'}_k + \mathds{1}(\bar{x}^{n'}_k = 1)(1 - x^{n'}_k)\right\} + \sum_{n' = n + 1}^{N}\sum_{k = 0}^{\bar{K}}\left\{\mathds{1}(\bar{z}^{n'}_k = 0)z^{n'}_k + \mathds{1}(\bar{z}^{n'}_k = 1)(1 - z^{n'}_k)\right\} \geq 1, \label{lbd}
\end{align}
where $\mathds{1}(x)$ is an indicator function that equals 1 if $x$ is true and 0 otherwise.

\subsection{Column Generation}
\label{sec:cg}
Problem \textbf{RF1} has an exponential number of variables, so it can be computationally prohibitive to enumerate all of them for reasonable-size instances. Therefore, we propose a column generation to solve it iteratively. First, problem \textbf{RF1} is initialized with a small subset of variables, called restricted master problem (RMP). Then, the RMP is solved by the simplex method, whereas a pricing problem is solved to generate new variables with negative reduced cost. These new variables are added to the RMP, and after that, the RMP is resolved. This process repeats until no variables with negative reduced cost are generated. The pricing problem with respect to problem \textbf{RF1} is as follows:
\begin{align}
&\min_{r \in \mathcal{R}^n}~c_r - \sum_{i \in \mathcal{I}^n}\alpha_{i,r}\nu_i -  \sum_{n' = n + 1}^{N}\beta_{n',r}\mu_{n'} - \mu_n. \label{pp}
\end{align}

The pricing problem belongs to the elementary shortest path problems with resource constraints (ESPPRCs) \citep{feillet2004exact,irnich2005shortest}, which are commonly solved by label-setting algorithms \citep{righini2008new}. The label-setting algorithms are a class of dynamic programming approaches that solve the ESPPRCs by state prorogation. In our case, states (or labels) represent partial routes from the depot to certain locations. By probably defining the states, the label-setting algorithms can enumerate all feasible routes, and hence guarantee to find an optimal route. Meanwhile, the algorithms can be speeded up by using special dominance rules to identify and discard redundant states. In summary, label-setting algorithms consist of three basic components: state definition, extension functions and dominance rules. Beside these three basic components, we also design and incorporate several important techniques to accelerate the label-setting algorithms, including the bounded bidirectional search \citep{righini2006symmetry}, ng-route relaxation \citep{martinelli2014efficient}, and label pruning techniques. The details of the label-setting algorithm for solving the pricing problem (\ref{pp}) are presented in Appendix \ref{sec:ls}.

\subsection{Route Enumeration Algorithm}
\label{sec:rea}
An intuitive method for solving problem \textbf{SF1} is a branch-and-price algorithm based on the column generation in Section \ref{sec:cg}. However, branch-and-price algorithms may converge slowly if branching decisions do not have strong impacts on the model. Therefore, we propose an iterative route enumeration algorithm to exactly solve problem \textbf{SF1}. The idea of this algorithm is similar to column generation. It first iteratively enumerates all feasible routes that possibly constitute optimal solutions of problem \textbf{SF1}, and then solves problem \textbf{SF1} with the enumerated routes directly by an MIP solver. The target routes for enumeration are given by Lemma \ref{lemma:re}:

\begin{lemma} \label{lemma:re} Given an upper bound  $ub$ of problem \textbf{SF1}, the optimal cost $\phi(\textbf{RF1})$ of problem \textbf{RF1}, a dual optimal solution $(\pmb{\mu}, \pmb{\nu})$ of problem \textbf{RF1} and a route $r \in \mathcal{R}^n$, $r$ cannot be in any optimal solutions if it satisfies:
\begin{align}
&c_r - \sum_{i \in \mathcal{I}^n}\alpha_{i,r}\nu_i -  \sum_{n' = n + 1}^{N}\beta_{n',r}\mu_{n'} - \mu_n  > ub - \phi(\textbf{RF1}).
\end{align}
\end{lemma}
Lemma \ref{lemma:re} states that if problem \textbf{RF1} is solved, a route with reduced cost larger than the gap $ub - \phi(\textbf{RF1})$ can not be in any optimal solution of problem \textbf{SF1}. 
Notably, our route enumeration algorithm does not require an upper bound as input, but iteratively generates good upper bounds. Let $ub$ be the best bound found by the algorithm, and $\delta$ be the gap used for route enumeration. At first,  $ub$ is initialized by positive infinity, and $\delta$ is initialized by a control parameter $StepSize$. If $\delta$ is too small and the enumerated routes cannot constitute a feasible solution or a better solution than $ub$, $\delta$ is increased by $StepSize$. Otherwise, the optimal solution obtained by solving problem \textbf{SF1} is used to update $ub$. Once $\delta$ is greater than or equal to $ub - \phi(\textbf{RF1})$, an optimal solution to problem \textbf{SF1} is found. 
Note that the route enumeration algorithm is also used to solve for the optimal single-period delivery cost in the offline estimation stage of AJRP, in which a large number of single-period problem instances have to be solved. We present the pseudocode and the detailed description of the algorithm in the Appendix \ref{sec:reap}.

\section{Computational Results and Discussion} \label{sec:sim}
In this section, we evaluate the performance of AJRP on both real-world and synthetic data sets. We first introduce the data set, simulation setup, and benchmark policies. Then we analyze the computational and on-time performance of AJRP and discuss managerial implications to the on-time delivery operations management.

\subsection{Data Sets}
The main data set is collected from our industry partner, whereas the synthetic instances are simulated to serve as additional test examples.
\subsubsection{Partner's Data}
The grocery chain store shared its order data set for on-demand meal boxes. The data set contains the following information of each placed delivery order: 1) order time: the time when the order is placed; 2) order quantity: the number of items (meal boxes) in the order; 3) time window: the delivery time target;
4) longitude and latitude: the customer location; 5) cutoff time:
the provider has set a sequence of evenly distributed cutoff times $\{t_1, t_2, \dots\}$ with $\Delta = t_{n+1} - t_n = 15$ minutes, and all orders placed within $[t_n, t_{n+1})$ are batched together and share the same delivery time window. We also acquired the travel distance data between customer locations (including the depot) from the Baidu Map API. Because the majority of orders were collected during the lunch peak hours, we focus on the time period from 10:00 am to 11:30 am, covering seven decision epochs (cut-off times).  We use two consecutive weeks of the order data as training set, and the orders in the following week make up the test set. To reflect different supply scenarios with varying fleet sizes, we consider the number of drivers as small (35 drivers in total), medium (40 drivers in total), and large (45 drivers in total) relative to demand. 

\subsubsection{Synthetic Data} To examine the scalability and generalizability of our algorithms, we perform additional computational studies on a set of synthetic instances. The instances are generated by varying the number of decision epochs ($N$) and demand generating process. Specifically, the potential customer locations are uniformly distributed on a plane of 10 km $\times$ 10 km, where the depot is located at the center. The distance between customer locations is calculated with Euclidean distance. The number of items ordered at each location follows a Poisson distribution with rate $\lambda=2$, and as a result, not all potential customers will place orders. We consider $N=\{10, 20\}$ and nonstationary demand arrival processes by setting the number of potential customer locations $I_n$ as a function of $n$.
Specifically, when $N=10$, $I_n = 25 + 5n$ for $n=1,\dots,5$ and $I_n=80 - 5n$ for $n=6,\dots,10$.  They are designed to mimic the practical scenario when demand peaks within the planning horizon (i.e., as in the partner's data). For $N=20$, the demand pattern of $N=10$ is repeated twice. Consequently, instances with $N=20$ replicate scenarios with two demand peaks (e.g., lunch and dinner hours) in the planning horizon. Figures that illustrate the temporal pattern of synthetic data are included in Appendix \ref{sec:distribution}. 100 random instances are generated for each configuration. 

The delivery speed is assumed to be $20$ km/hour and the on-site service time $s$ is set to be 5 minutes. The driver capacity is assumed to be $Q = 20$ (items) and the delivery duration limit is set to be $L_{\max} = 40$ minutes, as suggested by the partner. On the synthetic instances, the total number of drivers varies between small (54), medium (59), and large (64).

\subsection{Policy Implementation and Benchmarks}
During the offline estimation stage of AJRP, we generated 100 samples for each period in the evaluation of the single-period cost function $ \mathbb{E}_{\mathbb Q^{n}} \left[\mathcal H^s(K^{n},\mathcal{I}^{n},\BFq^{n})\right] $. We solve the single-period problems at different values of $K^n \in [\underline{K}^n, I_n]$, where $\underline{K}^n$ is the minimum required number of drivers in period $n$ (i.e., according to the capacity and delivery duration constraint). 
After the single-period problem solutions are collected across all samples, the single-period cost function is estimated using the simple sample average. Specifically, let $\mathcal{S}_n$ ($|\mathcal{S}_n|=100$) be the set of samples in period $n$, $\BFY_{s,k}$ be the optimal  single-period solution of sample $s \in \mathcal{S}_n$ with $k$ drivers, and $u(\BFY_{s,k})$ be the cost of solution $\BFY_{s,k}$. If there exists no feasible solution for sample $s \in \mathcal{S}_n$ given $k$ drivers, set $\BFY_{s,k} = \emptyset$. Let $\bar{\mathcal{S}}_{n,k} = \{s \in \mathcal{S}_n~|~ \BFY_{s,k} \neq \emptyset\}$. To account for the infeasible scenarios properly,  $\mathbb{E}_{\mathbb Q^n}\left[\mathcal H^s(k,\mathcal{I}^n,\BFq^n)  \right]$ $(n = 1, \ldots, N, k = 1,\ldots,K)$  is computed as $\sum_{s \in \bar{\mathcal{S}}_{n,k}}u(\BFY_{s,k})/|\bar{\mathcal{S}}_{n,k}|$ if $|\bar{\mathcal{S}}_{n,k}| \geq \varrho |\mathcal{S}_n|$ and $+\infty$ otherwise. The parameter $\varrho$ can be interpreted as a pruning parameter that controls the conservativeness of estimation. We set it to be $1/3$ in the experiments. The values of $\omega_{m}^{n'}(k)$ are computed in a similar way.



We compare our approach to two main benchmark dispatching and routing policies. The first benchmark policy replicates the current myopic policy used by the practitioner, and the second benchmark policy is adapted from a heuristic policy proposed in the literature:
\begin{enumerate}
	\item Simple myopic policy (current practice). As described in Section \ref{sec:background}, the company is using a simple myopic policy that disregards future order information in delivery planning. This policy dispatches and routes drivers in a way that only optimizes for the current batch of orders, i.e., by minimizing $ \sum_{i \in \mathcal I^n} u_i(\BFY^n)$. 
    Because this policy may not always be feasible when the fleet size is small (e.g., due to the capacity constraint), we follow a standard practice to introduce a set of third-party drivers of unlimited size and with extra labor cost. The labor cost of a third-party driver is proportional to his/her work time (total delivery time). Let $\BFY^n = (\BFY^{n,f}, \BFY^{n,p})$, where $\BFY^{n,f}$ and $\BFY^{n,p}$ correspond to the routing decision of the full-time drivers and the third-party drivers, respectively. Let $w(\BFY^{n,p})$ be the total work time of third-party drivers following $\BFY^{n,p}$. Then the simple myopic policy is derived by solving the following program:
\begin{align}
\min\ &     \sum_{i \in \mathcal I^n} u_i(\BFY^{n,f},\BFY^{n,p}) + \rho w(\BFY^{n,p})\label{bp2:obj}, \\
s.t. \ &(\BFY^{n,f},\BFY^{n,p})\in \mathscr{D}(\mathcal I^n, \BFq^n,\BFzeta^n) \label{bp2:con},
\end{align}
where $\rho$ is a weight parameter that reflects the additional labor cost. Without loss of generality, we set it to 10 so the operator has strong incentives to dispatch its own drivers and avoid calling third-party drivers.
	\item Adaptive myopic policy. The second benchmark policy is adapted from \cite{liu2020time}, where future order information is considered, but the driver dispatching and routing decisions are decoupled completely. Specifically, in each period, we first determine the number of dispatched drivers by solving a scheduling problem (after taking out the routing decision from (\ref{ajrp:obj}) - (\ref{ajrp:const2})):
	\begin{align}
	\min_{x_k^n\in \{0,1\}}\quad &  \sum_{k=0}^{\bar K^n} \mathcal{H}^s\left(k,\mathcal{I}^n,\mathbf q^n\right)x_k^n+\sum_{{n'}=n+1}^N \sum_{k=0}^{\bar K^{n'}}  \mathbb{E}_{\mathbb Q^{n'}}\left[\mathcal H^s(k,\mathcal{I}^{n'},\BFq^{n'})  \right] x_k^{n'}  \\
	\mbox{s.t.} \quad & \sum_{m=n}^{n'} \sum_{k=1}^{\bar K^m} x_k^m\omega_m^{n'}(k) \le \bar K^{n'} - \zeta^n_{n'} , \quad \forall n'=n,\cdots,N \\
	& \sum_{k=1}^{\bar K^{n'}} x_k^{n'} =1, \quad \forall n'=n,\cdots,N.
	\end{align}
	Then the single-period routing model is solved subject to the dispatching schedule constraint respecting the derived dispatching decision. Similar to AJRP, the adaptive myopic policy solves the dispatching problem in every period after collecting the new order information, i.e., the dispatching decision is updated in a rolling-horizon fashion. However, this policy is myopic in the routing part because it ignores the interactions between routing and future order arrivals, and the routing solution is derived independently from the dispatching decision. 
	It can be viewed as a combination of adaptive dispatching and myopic routing, which improves on the simple myopic policy to adjust the dispatching schedule according to future order arrivals.

\end{enumerate}

We also test another relevant heuristic policy that minimizes driver travel time to better balance driver capacity across different periods, of which the result is presented in Appendix \ref{sec:mintravel}. The algorithms were implemented in Java using callbacks of ILOG CPLEX 12.5.1. All of the experiments were conducted on a Dell personal computer with an Intel E5-1607 3.10 GHz CPU, 32 GB RAM, and Windows 7 operating system. To ensure all the policies are solved to optimality, the time limit is set to one hour per decision epoch. Note, however, this time limit is redundant for AJRP, as we show below that the solution time to AJRP is mostly within a few minutes.

\subsection{Computational Performance} \label{subsec:comp}
We report the solution time of the developed algorithms for the offline estimation and online optimization stages of AJRP. In the offline estimation stage, the single-period cost function $\mathcal H^s(k,\mathcal{I}^{n},\BFq^{n})$ has to be evaluated for a potentially large number of instances. 
On the synthetic instances, the average solution time for the single-period problem is 0.86 seconds (with a maximum of 188.27 seconds), which illustrates the promising computational performance of the proposed route enumeration algorithm. Furthermore, because the multiple traveling repairman problem (MTRP) can be treated a special case of $\mathcal H^s(k,\mathcal{I}^{n},\BFq^{n})$, we also evaluate the computational performance of our algorithm on three sets of public MTRP instances from the literature. Table \ref{table:singlePeriod} reports the number of instances tested, optimally solved, and the average solution time on each class of instances, compared with two state-of-the-art methods. The time limit of the route enumeration algorithm is set to one hour, while the time limits of the other two approaches are set to 2 hours. The results demonstrate that the proposed algorithm outperforms the existing methods in both solution time and quality, which bodes well for other on-time delivery problems built on MTRP. The detailed comparisons of these three approaches on the MTRP instances are presented in Tables \ref{table:set1}, \ref{table:set2} and \ref{table:set3} of Appendix \ref{sec:mtrp}. 

\begin{table}[!htp]
	\renewcommand{\arraystretch}{0.8}
\caption{Computational Results of the MTRP Instances}
\label{table:singlePeriod}
\begin{center}
\scalebox{0.8}{
	\renewcommand{\arraystretch}{1.5}
    \begin{tabular}{cccccccccc}
    \toprule
    \multirow{3}[0]{*}{Class} & \multirow{3}[0]{*}{\parbox{1.5cm}{~~~Total\\ Instances}} & \multicolumn{2}{c}{Route Enumeration Algorithm} & & \multicolumn{2}{c}{\citet{nucamendi2016mixed}} & & \multicolumn{2}{c}{\citet{muritiba2021branch}} \\
    \cline{3-4} \cline{6-7} \cline{9-10}
          &       & Instances  & Average Time & & Instances  & Average Time & & Instances  & Average Time \\
          &       & Tested/Solved & (In seconds) & & Tested/Solved & (In seconds) & & Tested/Solved & (In seconds) \\
    \midrule
    LQL  & 180   & 180/180 & \textbf{1.32}  & & 180/180 & 30.46  & & 180/180 & 31.58  \\
    E  & 12    & 12/12 & \textbf{127.10}  & & 9/9   & 507.86  & & 12/8  & 624.98  \\
    P  & 23    & 23/21 & \textbf{8.45}  & & 19/17 & 266.93  & & 23/17 & 1008.23  \\
    \bottomrule
    \end{tabular}%
}
\end{center}
\end{table}

Besides, we evaluate the solution efficiency of the proposed Benders decomposition framework for AJRP on the tested instances, and the results are summarized in Table \ref{table:solutiontime}. The average solution time per decision epoch of our framework is under 1 minute across different configurations (the maximum instance-specific solution time is 3 minutes), which marks a considerable improvement over the direct MILP formulation with CPLEX (the CPLEX solution time is well above 1 hour). Therefore, AJRP is practically feasible because optimal solutions can be returned during order preparation, which often takes more than 10 minutes. In general, instances with larger fleet sizes can be solved more efficiently because the corresponding delivery routes are shorter, and the pricing problems are easier to solve.

\begin{table}[htbp]
	\centering
	\caption{Average (Minimum and Maximum) Solution Time of AJRP per Decision Epoch (Seconds)}
	\begin{tabular}{ccc}
		\toprule
		Fleet Size &
		Partners' Data &  Synthetic Data \\
		\midrule
		Small    & 33.68 (18.47, 45.50)   & 3.36 (0.41, 72.61)\\
		Medium   & 24.01 (15.76, 32.04)   & 21.65 (1.87, 184.36)\\
		Large    & 15.03 (10.11, 28.31)   & 5.72 (2.26, 54.12)\\
		\bottomrule
	\end{tabular}%
	\label{table:solutiontime}%
\end{table}%


\subsection{Delivery Performance Improvement} \label{subsec:ontime}
We compare the on-time delivery performance of AJRP with the benchmark policies on real-world and synthetic instances. The delivery performance (cost) is measured by the sum of the delivery time of customer orders and the potential travel time of third-party drivers. For the chosen fleet sizes, the use of third-party drivers is very minimal, so the delivery performance mainly captures the order delivery time. We evaluate the relative performance improvement of AJRP over the simple myopic policy and the adaptive myopic policy by $(C^{\text{Myopic}} - C^{\text{AJRP}})/C^{\text{Myopic}}$, where $C^{\text{Myopic}}$ and $C^{\text{AJRP}}$ are the delivery cost of the myopic policy (static or adaptive) and AJRP, respectively. 

\begin{figure}[htbp]
	\centering
	\includegraphics[width=0.65\linewidth]{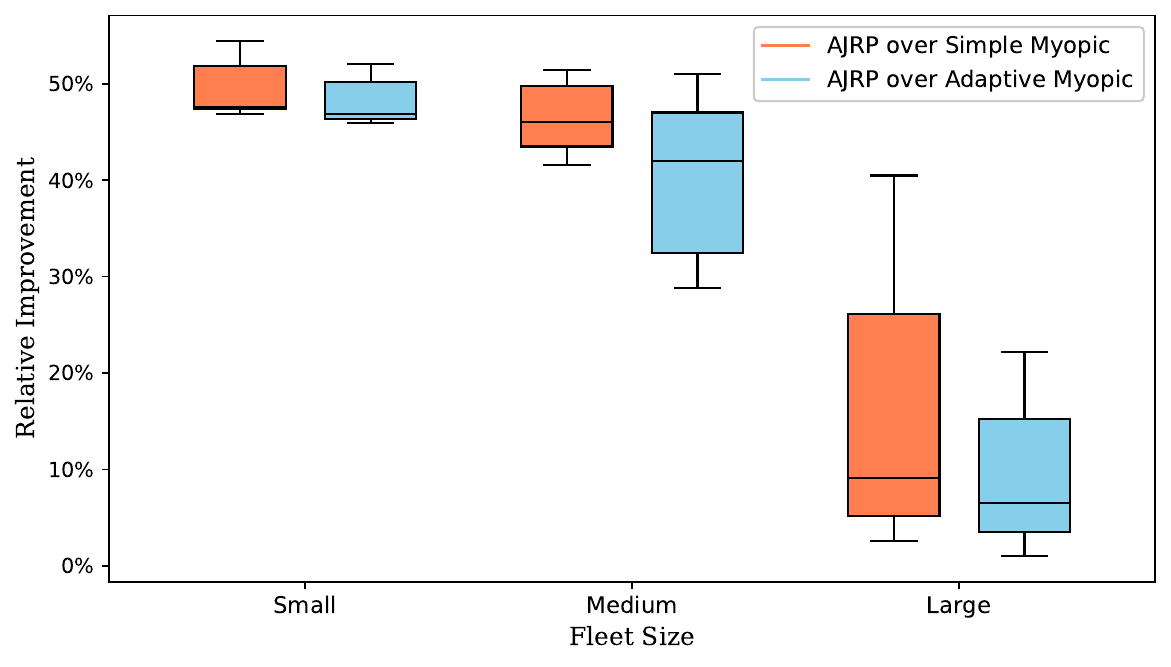}\hfill
	\caption{Relative Improvement of AJRP over Myopic Policies on Partner's Data}\label{fig:improve_real}
\end{figure}

Figure \ref{fig:improve_real} summarizes performance evaluation results of AJRP versus the two benchmark policies on partner's data. 
Compared to the current policy used by the company (simple myopic policy), AJRP provides an improvement of 36.53\% in delivery cost on average, which can translate to a substantial enhancement in delivery speed and promise reliability of on-demand orders. The average improvement of AJRP over the adaptive myopic policy is 32.29\%, which stresses the value of coordinating dispatching and routing decisions dynamically. Notably, these improvements are robust across different supply scenarios. Even when the driver supply is abundant, AJRP still significantly improves delivery performance.

\begin{figure}[htbp]
	\subfloat[$N=10$\label{sfig:improve_sim_n10}]{
		\includegraphics[width=0.5\linewidth]{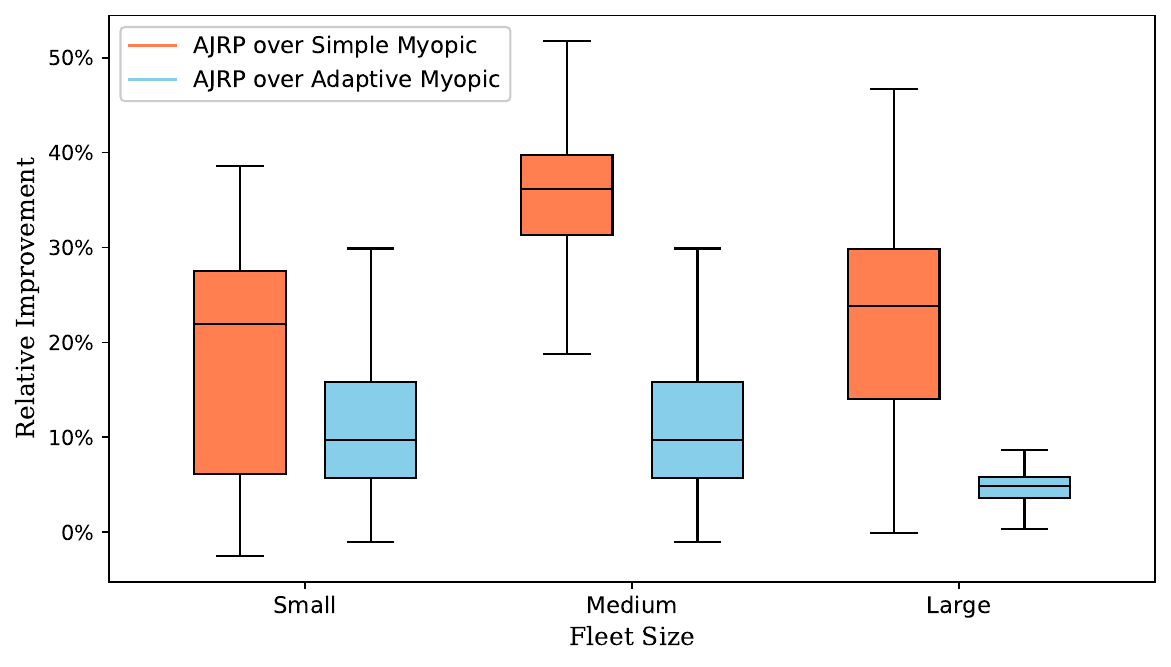}	}\hfill
	\subfloat[$N=20$\label{sfig:improve_sim_n20}]{
		\includegraphics[width=0.5\linewidth]{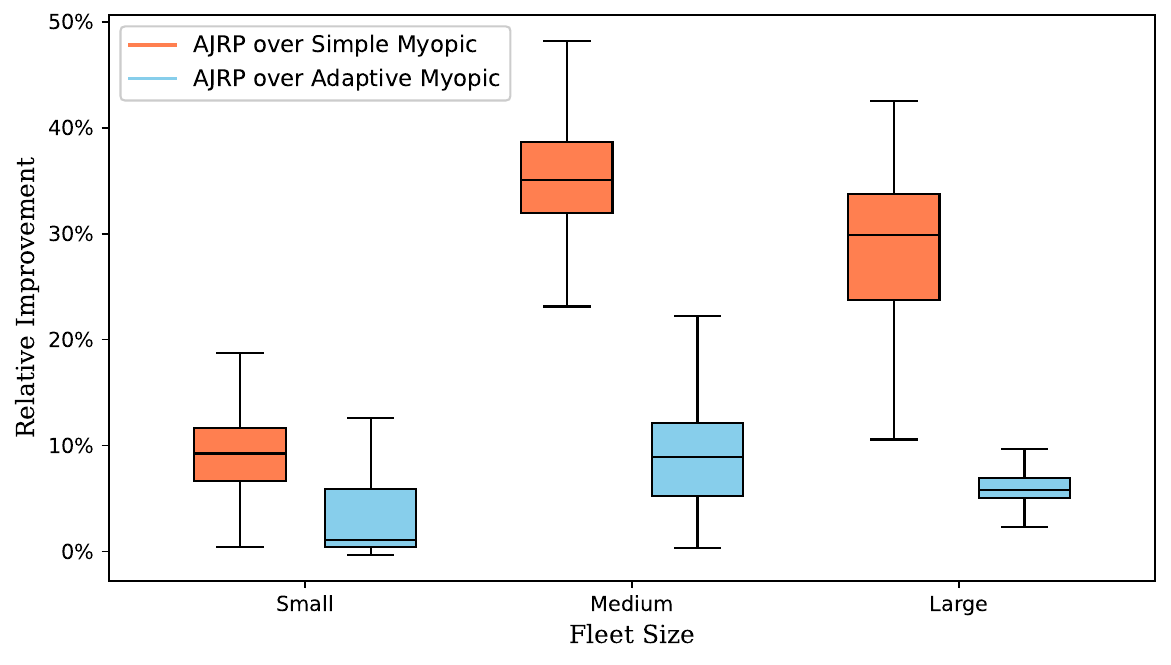}
	}\hfill
	\caption{Relative Improvement of AJRP over Myopic Policies on the Synthetic Data}\label{fig:improve_sim}
\end{figure}

Similar observations hold on the synthetic instances, of which the evaluation results are summarized in Figure \ref{fig:improve_sim}.   
Across different configurations, AJRP outperforms the two benchmark policies consistently. On average, AJRP outperforms the simple myopic policy by 24.82\% and the adaptive static policy by 8.63\% on the synthetic data. 
The improvement of AJRP tends to be greater for instances with medium fleet sizes than instances with small and large fleet sizes, in which dynamic optimization is more critical to matching supply and demand over time.

\subsection{Discussion and Policy Implications}
In this section, we perform several policy experiments and provide managerial insights for improving delivery performance based on partner's data.
\subsubsection{The Value of Dynamic Dispatching and Routing.}
Recall that the simple myopic policy follows both myopic dispatching and routing rules, whereas the adaptive myopic policy combines a dynamic dispatching rule with myopic routing. The improvement of the adaptive myopic policy over the simple myopic policy can be attributed to dynamic dispatching, and the improvement of AJRP over the adaptive myopic policy indicates the importance of dynamic routing. Therefore, we measure the relative value of dynamic dispatching and routing by the following two ratios: $(C^{\text{Simple Myopic}} - C^{\text{Adaptive Myopic}})/(C^{\text{Simple Myopic}} - C^{\text{AJRP}})$ and $(C^{\text{Adaptive Myopic}} - C^{\text{AJRP}})/(C^{\text{Simple Myopic}} - C^{\text{AJRP}})$, respectively. The higher the first ratio, the greater value dynamic dispatching generates (and the two ratios sum up to one). The average estimated values of these two ratios on the partner's data are presented in Table \ref{tab:value}. The main finding is that dynamic routing brings more benefits than dynamic dispatching for the company, and dynamic dispatching alone may not be sufficiently effective. However, as the fleet size gets larger, a higher contribution from dynamic dispatching can be observed, which implies dynamic dispatching is more valuable for large-fleet scenarios.

\begin{table}[htbp]
	\centering
	\caption{The Estimated Relative Value of Dynamic Dispatching and Routing}
	\begin{tabular}{ccc}
		\toprule
		Fleet Size & Value of Dynamic Dispatching & Value of Dynamic Routing \\
		\midrule
		Small    & 5.09\%   & 94.91\% \\
		Medium   & 23.70\%  & 76.30\% \\
		Large    & 40.02\%  & 59.98\% \\
		\bottomrule
	\end{tabular}%
	\label{tab:value}%
\end{table}%

\subsubsection{Comparing Dispatching and Routing Decisions.}
\begin{figure}[htbp]
	\centering
	\includegraphics[width=0.65\linewidth]{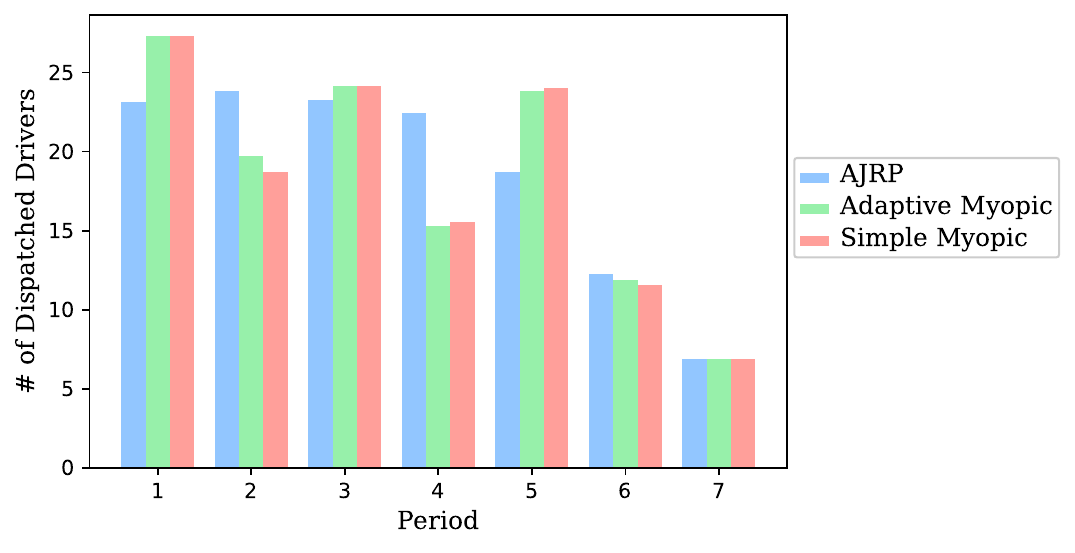}\hfill
	\caption{Dispatching Decision of Different Policies}\label{fig:driver_number}
\end{figure}

We first investigate the difference in the dispatching decision generated by the three policies. Figure \ref{fig:driver_number} presents the average number of dispatched drivers on the test set when the fleet size is 45 (large). We observe that the adaptive myopic policy behaves similarly to the simple myopic policy used by the company. This stresses that dynamic dispatching alone may not considerably impact the system performance. In contrast, AJRP dispatches drivers differently than the myopic policies: AJRP dispatches significantly fewer drivers in periods 1 and 5 but more drivers in periods 2 and 4. In particular, AJRP avoids sending out too many drivers in period 1 to better accommodate orders arriving in period 2. Although dispatching more drivers with shorter trips benefits on-time performance for the current batch of orders, blindly dispatching too many drivers poses risks of delaying future orders. This tradeoff is captured by AJRP more precisely than the myopic policies.

\begin{figure}[htbp]
	\subfloat[Period 1\label{sfig:duration_1}]{
		\includegraphics[width=0.45\linewidth]{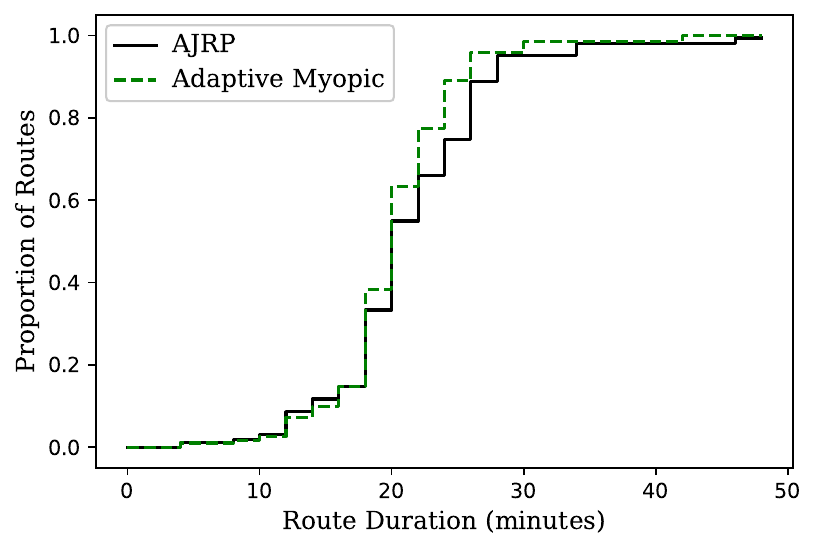}	}\hfill
	\subfloat[Period 3\label{sfig:duration_2}]{
		\includegraphics[width=0.45\linewidth]{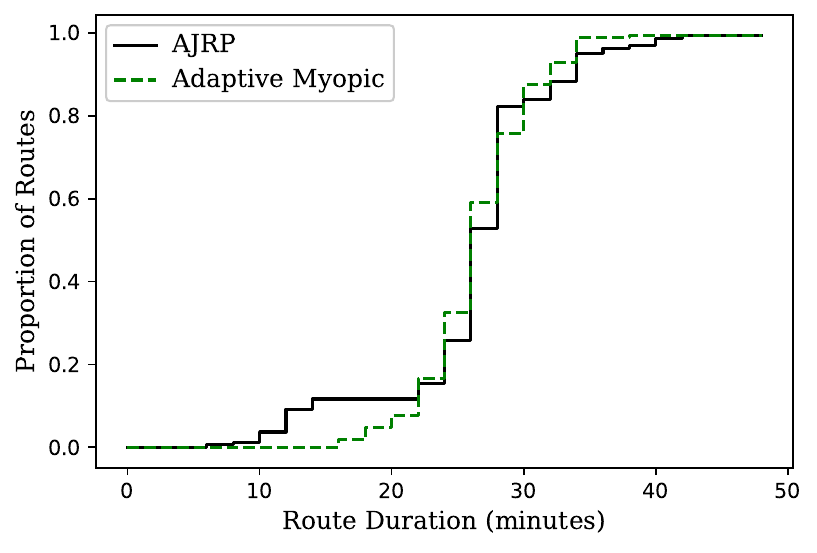}
	}\hfill
	\caption{Empirical Cumulative Distribution Function of Route Duration in Periods 1 and 3}\label{fig:route_duration_cdf}
\end{figure}

The interplay between dispatching decisions and future delivery performance lies in the planned delivery routes, of which the duration plays a major role in shaping future driver availability. In the considered setting, dispatched drivers whose routes are shorter than 15 minutes (30 minutes) can be dispatched again after one period (two periods). Figure \ref{fig:route_duration_cdf} depicts the empirical cumulative distribution function (CDF) of route duration under AJRP and the adaptive myopic policy for periods 1 and 3 (the simple myopic policy is omitted because it shares the same routing logic with the adaptive myopic policy). Note that AJRP dispatches fewer drivers in period 1, so we may expect longer routes from AJRP, and the dispatched drivers are less likely to return within the next two periods. However, due to careful routing optimization, the route duration of AJRP shares a similar distribution to that of the adaptive myopic policy. In particular, the percentage of routes that are shorter than 15 minutes and 30 minutes is almost the same under the two policies. We also compare the route duration in period 3, where the three policies dispatch a similar number of drivers. As shown in Figure \ref{sfig:duration_2}, AJRP plans more short routes (routes shorter than 15 minutes) than the adaptive myopic policy. Consequently, more drivers can be dispatched again in period 4 under AJRP, which boosts the overall on-time performance. These observations underline the value of routing optimization with multiple dispatch waves.

\subsubsection{Delivery Speed Versus the Fleet Size.} 
\begin{figure}[htbp]
	\centering
	\includegraphics[width=0.5\linewidth]{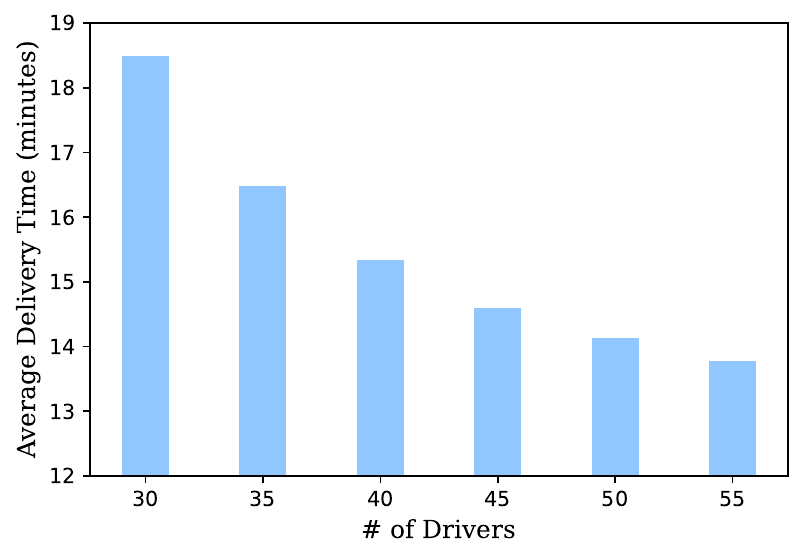}\hfill
	\caption{Average Delivery Time with Varying Fleet Sizes}\label{fig:delivery_time}
\end{figure}

As the on-demand delivery market becomes more competitive, the system operator can pursue faster deliveries with a larger fleet size. Figure \ref{fig:delivery_time} presents the evolution of average delivery time as a function of the number of drivers. If the company sets a 15-minute delivery time target (it corresponds to a customer waiting time of 30-35 minutes after accounting for order preparation and packaging), the fleet size should be at least 40. Further, our results imply diminishing returns on increasing the fleet size. As the fleet size grows from 30 to 35, the average delivery time can be reduced by 2 minutes. But when the fleet size is already large (e.g., 50), the incremental reduction of average delivery time is only half a minute. In theory, there is a physical limit to the average delivery time pertaining to the delivery region, travel speed, and service times. Approaching the lower limit can be economically unviable for many operators because of the resulting high labor cost.

\subsubsection{Setting the Right Frequency of Dispatch Waves.} 
\begin{figure}[htbp]
	\centering
	\includegraphics[width=0.5\linewidth]{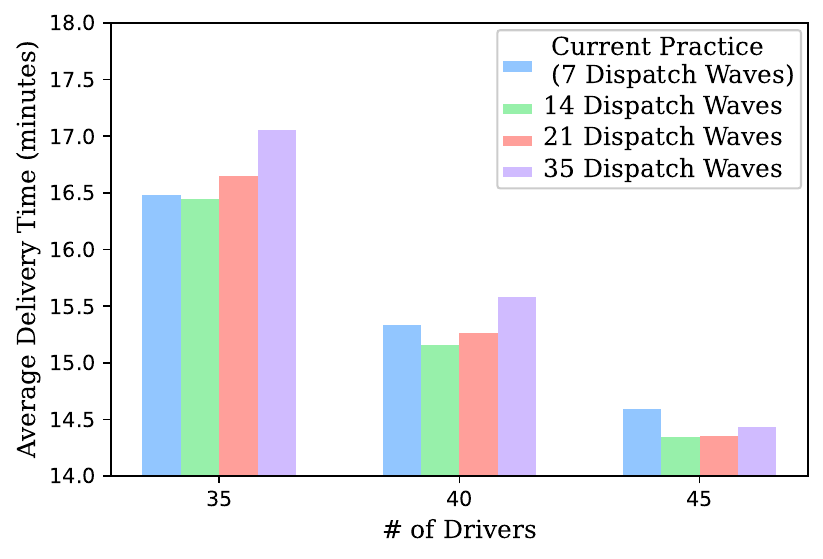}\hfill
	\caption{The Impact of Dispatch Frequency}\label{fig:dispatch_frequency}
\end{figure}

The company currently adopts 7 dispatch waves during peak hours. It is of high interest to the system operator to understand whether having more frequent dispatch waves (decision epochs) would benefit the on-time performance.  On the one hand, increasing the frequency of dispatch waves reduces the potential idle time of drivers and results in higher utilization of delivery capacity. On the other hand, making more frequent dispatches limits the potential of bundling orders (e.g., the opportunity that an order is bundled with a future order coming from a nearby location) and may compromise the route efficiency. To find the right frequency of dispatch waves, we increase the number of decision epochs to 14, 21, and 35, which implies a shorter time between consecutive epochs than the current practice. Figure \ref{fig:dispatch_frequency} depicts the average delivery time as a result of the increased dispatch frequency. The main observation is that using 14 dispatch waves obtains the best on-time performance: given a fleet of 45 drivers, the average delivery time can be reduced by 0.25 minutes when shifting from 7 dispatch waves to 14 dispatch waves. However, having too frequent dispatch waves may slow down the delivery process, as longer delivery time is observed for 21 and 35 dispatch waves.  Additionally, when the fleet size grows larger, the relative benefit from more frequent dispatches becomes more pronounced. This is because additional supply can be better utilized when the fleet is dispatched more frequently.

\subsubsection{The Value of Flexible Order Postponement.} \label{subsec:postpone} The considered dispatch policy, as motivated by the partner, does not allow flexible order postponement. Instead, an order is ``partially postponed" to be assigned at the next decision epoch and will not be postponed to later decision epochs. In theory, order postponement is beneficial when the postponed order can be effectively bundled with future orders and lead to more efficient delivery routes. The downside is that the postponed order has to wait for a longer time at the depot, which can negatively impact the on-time performance. Finding the optimal postponement strategy in a stochastic environment is challenging, so we examine the value of flexible order postponement in a clairvoyant manner. Assuming perfect future order information from the next period, we solve for the optimal order postponement decision along with driver dispatching and routing. Then we compare the average delivery time with and without order postponement (the dispatching decision follows AJRP when no order postponement is permitted). On the real instances with a medium number of drivers, the improvement in the average delivery time due to flexible order postponement varies from 1.41\% to 4.20\%, with an average of 2.56\%.\footnote{Based on our clairvoyant evaluation method, the estimated improvement from flexible order postponement is optimistic, and the actual improvement can be less than the reported value.} This suggests modest benefits from flexible order assignment to on-time performance, given that dispatching and routing decisions are well optimized. Nevertheless, jointly optimizing dispatching, routing, and postponement dynamically can be important for other metrics or applications, and we leave it as a future research direction.  

\subsubsection{The Impact of Sample Size.} The single-period cost function and the number of en route drivers are estimated from a simple sample average approximation. To understand how the sample size impacts the performance of AJRP, we vary the sample size $|\mathcal{S}_n|=10, 20,100,200$. We observe that the average delivery time does not vary significantly as the sample size changes. In particular, the average delivery time increases by up to 0.1 minutes when the sample size decreases from $100$ to $10$, with the largest increase observed for instances with medium fleet sizes. Interestingly, as the sample size grows from $100$ to $200$, the average delivery time does not improve, which can be attributed to overfitting the training sample.

\section{Concluding Remarks}
Fulfilling on-demand delivery orders rapidly in a dynamic and stochastic environment is a challenging task for many grocery and food retailers. The logistics system operator must dynamically optimize the dispatching and routing of drivers in response to new order arrivals and in anticipation of future orders. Computational difficulties are prominent due to the combinatorial nature of the problem and uncertain sequential arrivals of customer orders. Motivated by a large grocery chain store, we model and solve a stochastic dynamic dispatching and routing problem for on-time delivery of on-demand orders. We develop a structured approximation framework and computationally efficient algorithms that yield implementable solutions in real time. The proposed policy, AJRP, combines offline estimation and stochastic lookahead effectively. We show that AJRP enjoys a bounded approximation ratio and worst-case performance guarantee.
Our extensive computational experiments confirm the superior performance of AJRP on real-world and synthetic data sets. Compared to the current myopic policy used by the company, AJRP reduces the average delivery time by up to 49.61\%. 
Our results suggest dynamic routing is more beneficial than dynamic dispatching, especially when the fleet size is not so large. Due to the multi-trip and multi-dispatch features of our problem, a careful planning of routes plays an essential role in matching delivery capacity with demand. We also examine the impact of increasing dispatch frequency and the value of flexible order postponement. 

Our work has several limitations and can be extended in the following directions. First, because our modeling framework is focused on a single-depot setting where delivery orders originate from the central store, it would be interesting to consider a multi-depot scenario that allows order bundling across different stores. While our partner is not allowing such bundling policies because their stores are not in the vicinity of each other, some convenience stores may be able to explore the associated bundling flexibility. Second, it is possible to consider driver supply uncertainty in our model, which may be prominent when the company hinges on crowd-sourcing drivers to fulfill delivery orders. One may adjust the estimation procedure in our approximation framework accordingly -- the estimation of $\bar{\omega}_{m}^{n'}(K^m)$ can be tuned to reflect the case where crowd-sourcing drivers may not always return to the depot for future dispatch waves. Moreover, one can update the estimation of $\bar{\omega}_{m}^{n'}(K^m)$ by applying and evaluating the approximate dispatching policy iteratively to improve the framework. 
Lastly, integrating other decisions such as pricing and staffing with our model is interesting and may drive further methodological development.

\ACKNOWLEDGMENT{%
	The authors thank three anonymous referees, the associate editor, and Department Editor Melvyn Sim for their very timely and constructive comments. The authors acknowledge the support from the National Natural Science Foundation of China [Grants 72222011, 72171112],  the Young Elite Scientists Sponsorship Program by China Association for Science and Technology [Grant 2019QNRC001], and the Discovery Grant from the Natural Sciences and Engineering Research Council of Canada  [RGPIN-2022-04950].
}

\begingroup
	\setlength{\bibsep}{-0.1pt}
	\bibliographystyle{ormsv080}
	\bibliography{references}
\endgroup

\newpage

\section*{Online Appendix}
In this appendix, we present the key notations, technical proofs, additional algorithmic details, and more supporting evaluation results.
\begin{APPENDICES}
\section{Key Notations} \label{sec:notations}
\begin{table}[!htp]
	\renewcommand{\arraystretch}{1.2}
	\caption{Notations.}
	\begin{threeparttable}[t]
		\begin{center}
			\scalebox{0.75}{
				\begin{tabular}{cc}
					\toprule
					 Name & Description\\
					\midrule
				  \multicolumn{2}{c}{Driver Dispatching and Routing Model} \\
					\hline
					$N$ & Number of planning periods (decision epochs) \\
					$\Delta$ & Length of each planning period \\
					$t_n$ & Start time of period $n$ ($t_0$ refers to the start time of services)\\
					$t_p$ & Preparation time of orders (in a batch)\\
					$\bar{K}^n$ & Number of available drivers in period $n$ \\
					$\bar{K}$ & Total number of available drivers\\
					$v$ & Travel speed of drivers\\
					$Q$ & Vehicle capacity (maximum number of items loaded to a vehicle)\\
					$\mathcal{I}$ & Set of potential customer locations\\
					$\mathcal{I}^n$ & Set of customer locations realized between $[t_{n-1}, t_n)$\\
					$\BFq^n$ & Order quantity vector corresponding to 	$\mathcal{I}^n$\\	
					\multirow{2}{*}{$\BFzeta^n$} & Driver status vector: 
					$\BFzeta^n = (\zeta^n_1, \dots, \zeta^n_N)$, where $\zeta^n_{n'}$ is the number of en route drivers in period $n'$ \\
					& due to the dispatching decisions made prior to period $n$\\
					$\mathbb{Q}^{n}$ &  Joint distribution of the customer locations and order quantities in period $n$\\
					$\BFY^n$ & Joint dispatching and routing decision vector in period $n$\\
					$\mathscr{D}(\mathcal I^n, \BFq^n,\BFzeta^n)$ & Feasible region of $\BFY^n$\\
					$u_i(\cdot)$ & On-time performance measure for customer location $i$\\
					$L_{\max}$ & Hard delivery time target\\
					$l_k^n(\BFY^n)$ &  Route duration of driver $k$ dispatched in period $n$\\
					$\mathcal{H}_n(\mathcal I^n, \BFq^n,\BFzeta^n)$ & Cost-to-go function of period $n$ in the dynamic program\\
					\hline
					\multicolumn{2}{c}{Approximation Approach}\\
					\hline
					\multirow{2}{*}{$\mathcal H^s(K^n,\mathcal{I}^n,\mathbf q^n)$} & Single-period cost function with $K^n$ dispatched drivers when the realized customer locations and order\\
					&quantities are $\mathcal{I}^n$ and $\BFq^{n}$, respectively\\
					$\omega_{m}^{n'}(K^m)$ & Number of en route drivers in period $n'$ out of the $K^m$ drivers dispatched in period $m$\\
					\textbf{APT}$^{n}(\BFzeta^{n})$ & Approximation of  $\mathbb{E}_{\mathbb{Q}^{n}} \left[\mathcal{H}_{n}(\mathcal I^{n},  \BFq^{n}, \BFzeta^{n}) \right]$\\
					$x_k^{n}$ &  Binary decision variable that equals 1 if $k$ drivers are dispatched in period $n$ and 0 otherwise\\
					$\bar{\omega}_{m}^{n'}(K^m)$ & Estimated value of $\omega_m^{n'}(K^m)$\\
					$\bar{V}_{APT}^n(\BFzeta^{n})$ & Optimal objective value of \textbf{APT}$^n(\BFzeta^{n})$ with the choice of $\{\bar{\omega}_{m}^{n'}(K^m)\}_{\forall (m,n')}$\\
					\hline
					\multicolumn{2}{c}{Benders Decomposition}\\
					\hline
					$\mathcal{R}^n$ &  Set of feasible routes to serve orders in period $n$\\
					$c_r$ & Cost of route $r$\\
					$\alpha_{i,r}$ & Constant parameter that equals 1 if location $i$ is served in route $r$\\
					$\beta_{n',r}$ & Constant parameter that equals 1 if route $r $ is running in period $n'$\\
					$\theta_r$ & Binary decision variable that equals 1 if route $r$ is assigned to a driver and 0 otherwise\\
					\multirow{2}{*}{$z^{n'}_k$} & Binary decision variable that equals 1 if there exist $k$ drivers who are dispatched in period $n$ \\
					& but are still occupied in period $n'$\\
					$\mu_{n'}$ & Dual variables of constraints (25)\\
					$\mu_n$ & Dual variables of constraints (19)\\
					$\nu_i$ & Dual variables of constraints (26)\\
					$\pmb{\Theta}$  &  Set of extreme points of problem \textbf{SF1}'s dual problem\\
					$\pmb{\Lambda}$ & Set of extreme rays of problem \textbf{SF1}'s dual problem\\
					$\phi(\textbf{RF1})$ & Optimal cost of problem \textbf{RF1}\\
					\bottomrule
				\end{tabular}%
				\label{tab:notations}%
			}
		\end{center}
	\end{threeparttable}%
\end{table}

\section{Visualization of the Studied On-Demand Delivery Process} \label{sec:visual}
Figure \ref{fig:batch} illustrates a batch of orders and the corresponding delivery process in the studied on-demand delivery problem.
\begin{figure}[htbp]
	\centering
	\includegraphics[width=0.83\linewidth]{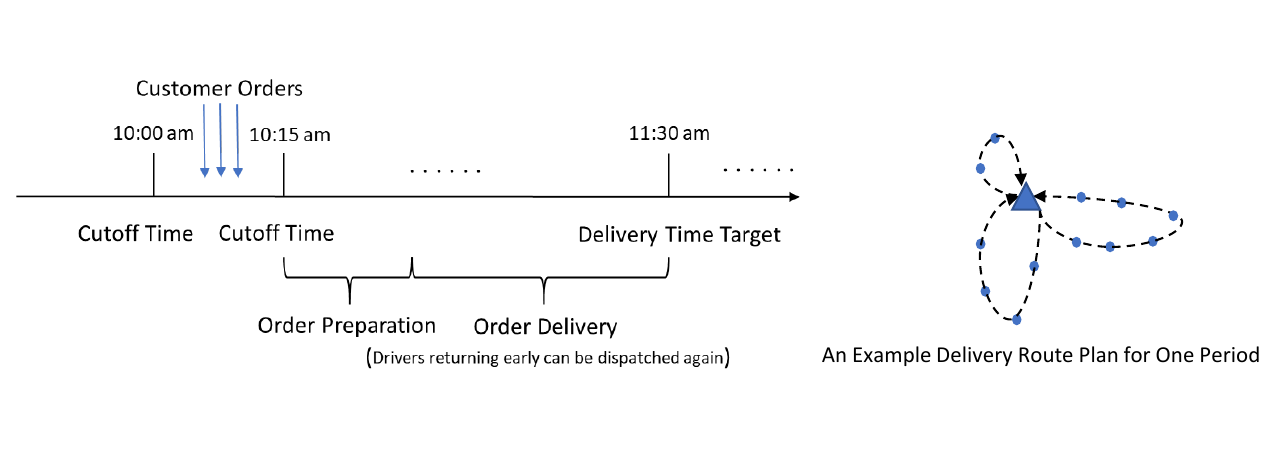}\hfill
	\caption{The Order Delivery Process: Triangle Represents the Store and Points Indicate Customer Locations}\label{fig:batch}
\end{figure}

\section{Constraints for $u_i(\BFY^n)$}
\label{sec:hdd}
Let  $t^n_{ij}$ be the travel time from location $i \in \mathcal{I}^n$ to location $j \in \mathcal{I}^n$, $s$ be the on-site service time, and $a^n_{ik}$ be a non-negative variable that records the delivery time of order $i \in \mathcal{I}^n$ by vehicle $k = \zeta^n_n + 1,\ldots, \bar{K}$. The following constraints are added for $u_i(\BFY^n)$ (the dependence on $\BFY^n$ is dropped below):
\begin{align}
a^n_{jk} \geq a^n_{ik} + s + t^n_{ij} + (y^n_{ijk} - 1)M, &\quad \forall i,j\in\mathcal{I}^n, k=\zeta^n_n + 1,\ldots, \bar{K}, \label{const:arrive} \\
a^n_{ik} \leq u_i \leq L_{\text{max}}, &\quad \forall i\in\mathcal{I}^n, k=\zeta^n_n + 1,\ldots, \bar{K}, \label{const:deadline}
\end{align}
where $M$ is a sufficiently large number. Constraints (\ref{const:arrive}) ensure that if vehicle $k$ delivers order $j$ immediately after order $i$, the delivery time of order $j$ should not be no earlier than the delivery time of order $i$ plus the travel time and the service time of order $i$. Constraints (\ref{const:deadline}) guarantee the delivery time of any order is no later than the target.

\section{Main Proofs}
\begin{proof}{Proof of Lemma \ref{lem:2}}
	(i) We prove the first part by induction. For the last period ($n=N$), $\mathcal{H}_N(\mathcal I^N, \BFq^N,\BFzeta^N)$ is nondecreasing in $\BFzeta^N$ because $\mathscr{D}(\mathcal I^N, \BFq^N,\BFzeta_1^N) \subseteq \mathscr{D}(\mathcal I^N, \BFq^N,\BFzeta_2^N)$ for $ \BFzeta_1^N \geq \BFzeta_2^N$. Now assume $\mathcal{H}_{n+1}(\mathcal I^{n+1}, \BFq^n,\BFzeta^{n+1})$ is nondecreasing in $\BFzeta^{n+1}$ for $n\leq N-1$. Then for any given $\mathbf{Y}^n \in \mathscr{D}(\mathcal I^n, \BFq^n, \BFzeta^n)$, $\BFzeta^{n+1}$ is nondecreasing in $\BFzeta^n$, and so $\mathcal{H}_{n+1}(\mathcal I^{n+1}, \BFq^{n+1},\BFzeta^{n+1})$ is nondecreasing in $\BFzeta^n$ by induction. It follows that the objective function of (\ref{eqn:bell1}) is nondecreasing in $\BFzeta^n$ (the expectation preserves the monotonicity). Therefore, $\mathcal{H}_{n}(\mathcal I^{n}, \BFq^{n},\BFzeta^{n})$ is nondecreasing in $\BFzeta^{n}$ by the property of $\mathscr{D}(\mathcal I^n, \BFq^n,\BFzeta^n)$, finishing the inductive proof.
	
	(ii) The second part follows from the fact that the feasible region of \textbf{APT}$^n(\BFzeta^n)$ is nonincreasing in $\BFzeta^n$.  \halmos
\end{proof}

\begin{proof}{Proof of Theorem \ref{them:approx}}
	(1) \textbf{Lower bound.} We can set $\omega_{n'}^{n'}(K^{n'}) = K^{n'}$ for $n' = n,\dots, N$ and $\omega_{m}^{n'}(K^{m}) = 0$ for $m<n'$. This implies that drivers will be able to return to the depot within one period of their dispatch. Therefore, under this set of $\omega_{m}^{n'}(K^m)$, we have a system wherein all the drivers are available at the beginning of every period, which is clearly a lower bounding system. More specifically, the optimal solution to \textbf{APT}$^n(\BFzeta^{n})$ would satisfy $K^{n'} = \bar{K} - \zeta^n_{n'}$, of which the objective value $\sum_{n'=n}^N   \mathbb{E}_{\mathbb Q^{n'}} \left[\mathcal H^s(\bar{K} - \zeta^n_{n'},\mathcal{I}^{n'},\BFq^{n'})\right]$ will be less or equal to $\mathbb{E}_{\mathbb{Q}^{n}} \left[\mathcal{H}_{n}(\mathcal I^{n},  \BFq^{n}, \BFzeta^{n}) \right]$ by the definition of $\mathcal{H}^s(\cdot,\cdot,\cdot)$.
	
	(2) \textbf{Upper bound.} By the existence assumption, we can find a feasible myopic policy $\pi^{sm}$ that dispatches $K^{n', sm}$ drivers in period $n'$ regardless of the realization of $(\mathcal I^{n'},  \BFq^{n'}, \BFzeta^{n'})$.  Then we can choose the values of $\omega_{m}^{n'}(K^m)$ such that the optimal solution to \textbf{APT}$^n(\BFzeta^{n})$ is $\{K^{n', sm}\}_{n' = n,\dots, N}$, and the optimal objective value is the same as the expected cost of $\pi^{sm}$ (one such choice is to set $\omega_{m}^{n'}(K^{m,sm}) = 0$ and $\omega_{m}^{n'}(K) = \bar{K} + 1$ for $K\neq K^{m,sm}$). Because $\pi^{sm}$ is a feasible policy and its expected cost will be an upper bound of $\mathbb{E}_{\mathbb{Q}^{n}} \left[\mathcal{H}_{n}(\mathcal I^{n},  \BFq^{n}, \BFzeta^{n}) \right]$.
	
	Lastly, when $n =N$, we have
	\begin{align*}
	\mathbb{E}_{\mathbb{Q}^{n}} \left[\mathcal{H}_{n}(\mathcal I^{n},  \BFq^{n}, \BFzeta^{n}) \right] = \mathbb{E}_{\mathbb{Q}^{N}} \left[\mathcal{H}_{N}(\mathcal I^{n},  \BFq^{N}, \BFzeta^{N}) \right] = \min_{K^N \in \mathbb{N}} & \mathbb{E}_{\mathbb{Q}^{N}} \left[\mathcal{H}^{s}(K^N, \mathcal I^{n},  \BFq^{N}) \right] \\
	s.t. \ & K^N \leq \bar{K} - \zeta_N^N,
	\end{align*}
	which is the same as problem \textbf{APT}$^n(\BFzeta^{n})$ because $\omega_N^N(K^N) = K^N$. \halmos
\end{proof}

\begin{proof}{Proof of Corollary \ref{lem:n2}}
	Due to Theorem \ref{them:approx}, the approximation used in AJRP is exact for the last period. Therefore, AJRP is optimal for the two-period JDR. \halmos
\end{proof}
\begin{proof}{Proof of Proposition \ref{prop:approx}}
	When there exists a static myopic policy, \textbf{APT}$^n(\BFzeta^{n})$ is feasible with $\omega_m^{n'}(K^m) = \bar{\omega}_{m}^{n'}(K^m)$ so $\bar{V}_{APT}^n(\BFzeta^{n})$ is finite. Let define
	\begin{align}
	\bar{V}_{APT}^n(\BFzeta^{n}, \mathbf{\epsilon}) := \min_{K^{n'} \in \mathbb{N}} & \sum_{n'=n}^N   \mathbb{E}_{\mathbb Q^{n'}} \left[\mathcal H^s(K^{n'},\mathcal{I}^{n'},\BFq^{n'})\right] \nonumber \\
	s.t.\ & \sum_{m=n}^{n'} \bar{\omega}_{m}^{n'}(K^m)  \leq \bar{K} - \zeta^n_{n'} + \epsilon(n'),\quad \forall n'=n,\dots, N. \label{approxprob2}
	\end{align}
	We now prove that there exist $\mathbf{\epsilon}_1, \mathbf{\epsilon}_2\in \mathbb{R}^{N-n+1}$ satisfying $\mathbf{\epsilon}_1\geq 0$ and  $\mathbf{\epsilon}_2\leq 0$, such that
	\begin{align*}
	\bar{V}_{APT}^n(\BFzeta^{n}, \mathbf{\epsilon}_1) \leq \mathbb{E}_{\mathbb{Q}^{n}} \left[\mathcal{H}_{n}(\mathcal I^{n},  \BFq^{n}, \BFzeta^{n}) \right] \leq \bar{V}_{APT}^n(\BFzeta^{n}, \mathbf{\epsilon}_2).
	\end{align*}
	\begin{enumerate}
		\item $\bar{V}_{APT}^n(\BFzeta^{n}, \mathbf{\epsilon}_1) \leq \mathbb{E}_{\mathbb{Q}^{n}} \left[\mathcal{H}_{n}(\mathcal I^{n},  \BFq^{n}, \BFzeta^{n}) \right]$: we can choose a large enough $\epsilon_1\geq 0$ such that it is possible to dispatch $I$ drivers in every period (recall $I$ is the maximum number of customer locations), which yields the best possible on-time performance and thus is less than $\mathbb{E}_{\mathbb{Q}^{n}} \left[\mathcal{H}_{n}(\mathcal I^{n},  \BFq^{n}, \BFzeta^{n}) \right]$.
		\item  $\mathbb{E}_{\mathbb{Q}^{n}} \left[\mathcal{H}_{n}(\mathcal I^{n},  \BFq^{n}, \BFzeta^{n}) \right] \leq \bar{V}_{APT}^n(\BFzeta^{n}, \mathbf{\epsilon}_2)$: for period $n'=n,\dots, N$, there exists a lower bound $\underline{K}^{n'}$ such that $\mathbb{E}_{\mathbb Q^{n'}} \left[\mathcal H^s(K^{n'},\mathcal{I}^{n'},\BFq^{n'})\right]$ is only finite when $K^{n'} \geq \underline{K}^{n'}$. Now we construct $\mathbf{\epsilon}_2$ by setting $\epsilon_2(n')= \sum_{m=n}^{n'}   \bar{\omega}_{m}^{n'}(\underline{K}^m) - (\bar{K} - \zeta^n_{n'}),\  \forall n'=n,\dots, N $. By doing so the only feasible solution to problem (\ref{approxprob2}) is $\{\underline{K}^{n'}\}_{n' = n,\dots, N}$ (from the property that $\bar{\omega}_{n'}^{n'}(\cdot)$ is a non-decreasing function). As a result, $\bar{V}_{APT}^n(\BFzeta^{n}, \mathbf{\epsilon}_2) = \sum_{n'=n}^N   \mathbb{E}_{\mathbb Q^{n'}} \left[\mathcal H^s(\underline{K}^{n'},\mathcal{I}^{n'},\BFq^{n'})\right] $, which is greater than the cost of any feasible static myopic policy $\pi^{sm}$ because $K^{n',sm}\geq \underline{K}^{n'}$. Consequently, $\bar{V}_{APT}^n(\BFzeta^{n}, \mathbf{\epsilon}_2) \geq  \mathbb{E}_{\mathbb{Q}^{n}} \left[\mathcal{H}_{n}(\mathcal I^{n},  \BFq^{n}, \BFzeta^{n}) \right] $.
	\end{enumerate}

	Next, let $\vartheta = \bar{V}_{APT}^n(\BFzeta^{n}, \mathbf{\epsilon}_2)/\bar{V}_{APT}^n(\BFzeta^{n}, \mathbf{\epsilon}_1)$, then we have
	\[ \frac{1}{\vartheta} \leq  \frac{\bar{V}_{APT}^n(\BFzeta^{n})}{\mathbb{E}_{\mathbb{Q}^{n}} \left[\mathcal{H}_{n}(\mathcal I^{n},  \BFq^{n}, \BFzeta^{n}) \right]} \leq  \vartheta\]
	because  $\bar{V}_{APT}^n(\BFzeta^{n}, \mathbf{\epsilon}_1) \leq \bar{V}_{APT}^n(\BFzeta^{n}) \leq \bar{V}_{APT}^n(\BFzeta^{n}, \mathbf{\epsilon}_2)$.
	Lastly, when $\bar{K}$ is large enough for $n'=n,\dots, N$, the optimal dispatching and routing policy would be setting $K^{n'} = I$ in every state, i.e., each driver will serve at most one customer location, which corresponds to the best achievable static myopic policy. As such $\bar{V}_{APT}^n(\BFzeta^{n}) = \mathbb{E}_{\mathbb{Q}^{n}} \left[\mathcal{H}_{n}(\mathcal I^{n},  \BFq^{n}, \BFzeta^{n}) \right]$.\halmos

\end{proof}

\begin{proof}{Proof of Theorem \ref{prop:bound}}
	We first focus on a single period with $I^*$ customer locations. As done in the proof of Proposition \ref{prop:approx}, the shortest expected delivery time is achieved when we can dispatch $I^*$ drivers so each driver only serves one customer. As such, the expected delivery time is $\bar{r}/v+s$ per order. This is the lower bound of the single-period expected delivery time for $I^*$ customer orders, and we denote it by $V_{lb}(I^*)$.
	
	Due to the capacity constraint (and the possible delivery deadline), the number of dispatched drivers is at least $l= \lceil I^*/Q \rceil$. Now we provide an upper bound for the expected delivery time when dispatching $l$ drivers. The routing policy follows the well known tour partitioning scheme proposed by \cite{haimovich1985bounds}. Under this scheme, we first construct the optimal tour (TSP tour, $\mathcal{TSP}(I^*)$) through all the customer locations. Then we split the tour into $l$ segments and create $l$ feasible driver tours by connecting the depot with the endpoints of the segments.
	
	Conditioning on the locations of the $I^*$ customers, there are $I^*$ ways to partition the optimal tour, each corresponding to a different starting location. Among the $I^*$ different partitions, a customer location will be connected to the depot (as the first visit location) $l$ times. When location $i$ is selected as the first visit location, the radial travel time will contribute at most $Q\cdot r_i/v$ to the delivery time of orders on the tour starting at $i$ (recall that a driver tour contains $Q$ orders). Summing over the $I^*$ possible partitions, the total contribution from radial travel time is $l\cdot Q\sum_{i=1}^{I^*} r_i /v $. Next we compute the contribution from the TSP travel time. For a given arc on the tour, $(i,i+1)$, its contribution to the total delivery time depends on the delivery sequence of $i$ on the driver tour. When $i$ is the $k$th visited customer on the tour, arc $(i,i+1)$ will contribute $(Q-k)r(i,i+1)/v$ to the total delivery time (the customers following $i$ all include $r(i,i+1)/v$ as part of their delivery time). Because $k$ can take values from $1$ to $Q$ ($Q$ indicates $(i,i+1)$ is not on any driver tours, i.e., $i$ and $i+1$ are connecting to the depot), the total contribution from the TSP travel time over all the possible partitions is
	\[  l \sum_{(i,i+1)\in \mathcal{\mathcal{TSP}(I^*)}}  \sum_{k=1}^{Q} (Q-k) r(i,i+1)  = \frac{l\cdot Q(Q-1)}{2v} L(\mathcal{\mathcal{TSP}(I^*)}), \]
	where $L(\mathcal{TSP}(I^*))$ is the TSP tour length through the $I^*$ locations. Following a similar argument, the total contribution from on-site service time is $l\cdot Q(Q+1)s/2$ for each partition. Therefore, the total delivery time of customer orders summing over the $I^*$ different partitions is
	\begin{align*}
	\frac{l\cdot Q\sum_{i=1}^{I^*} r_i }{v} + \frac{l\cdot Q(Q-1)}{2v} L(\mathcal{\mathcal{TSP}(I^*)}) + \frac{I^*\cdot l\cdot   Q(Q+1)s}{2}.
	\end{align*}
	With a randomly selected partitioning, the expected delivery time would be
	\begin{align*}
	\frac{\sum_{i=1}^{I^*} r_i }{v} + \frac{(Q-1)}{2v} L(\mathcal{\mathcal{TSP}(I^*)}) + \frac{I^*(Q+1)s}{2},
	\end{align*}
	where we utilize the relationship $l\cdot Q\approx I^*$. Then we take expectation with respect to the random locations of customers and derive the expected delivery time as
	\begin{align*}
	\frac{I^* \bar{r} }{v} + \frac{(Q-1)}{2v} \mathbb{E}(L(\mathcal{\mathcal{TSP}(I^*)})) + \frac{I^*(Q+1)s}{2}.
	\end{align*}
	
	Because the above tour partitioning policy does not necessarily minimize the expected delivery time and only obtains a feasible solution, it provides an upper bound of the expected delivery time for a period with $I^*$ customer locations,  $V_{ub}(I^*)$. Then we have
	\begin{align*}
	\frac{V_{ub}(I^*)}{V_{lb}(I^*)} = \frac{I^* \bar{r}/v  + (Q-1)  \mathbb{E}(L(\mathcal{\mathcal{TSP}(I^*)})/2v + I^*(Q+1)s/2}{(\bar{r}/v+s)I^*}.
	\end{align*}
	The Beardwood-Halton-Hammersley (BHH) theorem implies that for large $I^*$
	\[  \frac{\mathbb{E}(L(\mathcal{\mathcal{TSP}(I^*)})}{I^*} \approx \beta \frac{\sqrt{A}}{\sqrt{I^*}} ,   \]
	where $\beta$ is the TSP constant \citep{beardwood1959shortest}. Hence, we can obtain an approximation for $V_{ub}(I^*)/V_{lb}(I^*)$ as
	\begin{align*}
	\frac{V_{ub}(I^*)}{V_{lb}(I^*)} \approx \frac{ \bar{r}/v  + \beta (Q-1) \sqrt{A}/(2v\sqrt{I^*}) + (Q+1)s/2}{\bar{r}/v+s}.
	\end{align*}
	
	Now we consider the original expression of $\vartheta$, which follows
	\[\vartheta = \frac{\bar{V}_{APT}^n(\BFzeta^{n}, \mathbf{\epsilon}_2)}{\bar{V}_{APT}^n(\BFzeta^{n}, \mathbf{\epsilon}_1)} = \frac{ \sum_{n'=n}^N   \mathbb{E} \left(\mathcal H^s(\underline{K}^{n'},I^{n'})\right)}{ \sum_{n'=n}^N \mathbb{E} (V_{lb}(I^{n'})) } \leq \frac{ \sum_{n'=n}^N   \mathbb{E}(V_{ub}(I^{n'}))}{ \sum_{n'=n}^N \mathbb{E} (V_{lb}(I^{n'})) },\]
	where the expectation is taken with respect to the number of customer locations. Because we assume each customer orders exactly one item, the order quantity information is redundant and removed from the single-period optimal delivery time function. Moreover, observing that $V_{ub}(I)/V_{lb}(I)$ is decreasing in $I$, we have
	\begin{align*}
	\frac{ \mathbb{E}(V_{ub}(I^{n'}))}{\mathbb{E} (V_{lb}(I^{n'})) } \leq \frac{V_{ub}(I^*)}{V_{lb}(I^*)}\quad \forall n'=n,\dots, N
	\end{align*}
	due to the assumption that $I^{n'} \geq I^*$ for any $n'$. It follows that
	\begin{align*}
	\vartheta \leq \frac{ \sum_{n'=n}^N   \mathbb{E}(V_{ub}(I^{n'}))}{ \sum_{n'=n}^N \mathbb{E} (V_{lb}(I^{n'}))} \leq \frac{V_{ub}(I^*)}{V_{lb}(I^*)} \approx \frac{ \bar{r}/v  + \beta (Q-1) \sqrt{A}/(2v\sqrt{I^*}) + (Q+1)s/2}{\bar{r}/v+s},
	\end{align*}
	which holds for large $I^*$. \halmos
\end{proof}

\section{Logic Benders Decomposition}
\label{sec:dbd}

The pseudocode of the logic Benders decomposition is presented in Algorithm \ref{alg:bd}. In this algorithm, $\epsilon$ is the tolerance of the optimality gap, and $lb$ ($ub$) is the best lower bound (upper bound) found by the algorithm. Let $\phi(\pmb{SF1})$, $\phi(\pmb{DF1})$ and $\phi(\pmb{F2})$ be the optimal costs of problems \textbf{SF1}, \textbf{DF1} and \textbf{F2}, respectively. In each iteration, the algorithm first solves problem \textbf{F2} and then solves problem \textbf{RF1}. Depending on the results, different types of Benders cuts are added to problem \textbf{F2} to cut off the incumbent solution. If \textbf{RF1} is infeasible, it indicates that the feasibility cut (\ref{f2:con2}) is able to cut off the incumbent solution. If $\phi(\pmb{RF1}) - \bar{\eta}$ is greater than or equal to $ub - lb$, it indicates that the optimality cut (\ref{f2:con1}) is able to cut off the incumbent solution. Otherwise, the logic Benders cut (\ref{lbd}) is used, and problem \textbf{SF1} is solved exactly  to update the best upper bound. The algorithm continues until the optimality gap is smaller than or equal to the given threshold.
\begin{algorithm}[!h]
\caption{{\em Benders Decomposition}} \label{alg:bd}
\begin{small}
\begin{algorithmic}[1]
\STATE $\pmb{\Theta} \leftarrow \emptyset, \pmb{\Lambda} \leftarrow \emptyset, lb \leftarrow 0, ub \leftarrow \infty$; \label{bd:1}
\WHILE{$ub - lb > \epsilon$}
  \STATE Solve problem \textbf{F2} and obtain the optimal solution $(\bar{\pmb{x}},\bar{\pmb{z}}, \bar{\eta})$ ; \label{bd:2}
  \STATE $lb \leftarrow \phi(\pmb{F2})$;  \label{bd:3}
  \STATE Construct problem \textbf{RF1} according to $(\bar{\pmb{x}},\bar{\pmb{z}})$; \label{bd:4}
  \STATE Solve problem \textbf{RF1} by column generation and obtain the dual optimal solution $(\pmb{\mu}, \pmb{\nu})$; \label{bd:5}
  \IF{problem \textbf{RF1} is infeasible}
    \STATE Add the feasibility cut (\ref{f2:con2}) with respect to $(\pmb{\mu}, \pmb{\nu})$ to problem \textbf{F2}; \label{bd:6}
  \ELSIF{$\phi(\pmb{RF1}) - \bar{\eta} \geq ub - lb$}
    \STATE Add the optimality cut (\ref{f2:con1}) with respect to $(\pmb{\mu}, \pmb{\nu})$ to problem \textbf{F2};\label{bd:7}
  \ELSE
  \STATE Add the logic Benders cut (\ref{lbd}) with respect to $(\bar{\pmb{x}},\bar{\pmb{z}})$ to problem \textbf{F2}: \label{bd:11}
    \STATE Construct problem \textbf{SF1} according to $(\bar{\pmb{x}},\bar{\pmb{z}})$; \label{bd:8}
    \STATE Solve problem \textbf{SF1} by a route enumeration based algorithm; \label{bd:9}
    \STATE $ub \leftarrow  \min\{ub, \phi(\pmb{F2}) - \bar{\eta} + \phi(\pmb{SF1})\}$; \label{bd:10}
  \ENDIF
\ENDWHILE
\end{algorithmic}
\end{small}
\end{algorithm}

\section{Label-Setting Algorithm}
\label{sec:ls}
\subsection{Basic Components}
Because the reduced cost of a route is influenced by the time periods it goes through, the routes can be divided into different classes according to the period when a route returns back to the depot. Each class of routes correspond to a pricing problem. Now consider the pricing problem with respect to routes which must return back to the depot no later than period $n'$. Let $T_{\text{max}} = t_{n'} - t_n$ be the maximum duration of the routes. Each pair of locations $i$ and $j$ $(i,j \in \mathcal{I}^n \cup \{0\}, i \neq j)$ is associated with a cost $\bar{c}^n_{ij}$ as follows:
\begin{align}
&\bar{c}^n_{ij} =
\begin{cases}
-\frac{\nu_i}{2} - \frac{\nu_j}{2}, &~\text{if}~i,j \neq 0\\
-\frac{\sum_{k = n}^{n'}\mu_k}{2} - \frac{\nu_j}{2}, &~\text{if}~i = 0, j \neq 0\\
-\frac{\nu_i}{2} -\frac{\sum_{k = n}^{n'}\mu_k}{2}, &~\text{if}~ i \neq 0, j = 0.
\end{cases}
\end{align}

\textit{Label Definition}. Let $L_i = (\bar{c}(L_i), e(L_i), d(L_i), \mathcal{V}(L_i))$ be a label representing a partial path from depot to location $i \in \mathcal{I}^n \cup \{0\}$ where
\begin{itemize}
  \item $\bar{c}(L_i)$ is the reduced cost of the path;
  \item $e(L_i)$ is the earliest arrival time at location $i$;
  \item $d(L_i)$ is the total demand of the visited locations;
  \item $\mathcal{V}(L_i)$ is the set of locations that the path can extend to.
\end{itemize}

\textit{Extension functions.} The extension starts with an initial label $L_0 = \{0, 0, 0, \mathcal{I}^n \cup \{0\}\}$. For a pair of locations $i$ and $j$ and a label $L_i$ associated with location $i$, label $L_i$ can be extended to location $j$ to create a new label $L_j$ by the following extension functions:
\begin{align}
&e(L_j) = e(L_i) + s^n_i + t^n_{ij}, \label{ext:1}\\
&\bar{c}(L_j) =
\begin{cases}
\bar{c}(L_i) + e(L_j) + \bar{c}^n_{ij}, &\quad \text{if}~j \neq 0\\
\bar{c}(L_i) + \bar{c}^n_{ij}, &\quad \text{if}~j = 0
\end{cases} \label{ext:2}\\
&d(L_j) = d(L_i) + q^n_j, \label{ext:3}\\
&\mathcal{V}(L_j) = \left\{k \in \mathcal{V}(L_i) \backslash \{j\}~|~d(L_j) + q^n_k \leq Q, e(L_j) + s^n_j + t^n_{jk} \leq \min\{L_{\text{max}}, T_{\text{max}} - s^n_k - t^n_{k0}\}\right\}, \label{ext:4}
\end{align}
where $q^n_j$ is the demand of order $j$ and $Q$ is the vehicle capacity.

\textit{Dominance rules.} Let $L^1_i$ and $L^2_i$ be two labels associated with location $i$. Then $L^1_i$ dominates $L^2_i$ if they satisfy the following conditions:
\begin{align}
&\bar{c}(L^1_i) \leq \bar{c}(L^2_i), \label{dom:1}\\
&e(L^1_i) \leq e(L^2_i), \label{dom:2}\\
&d(L^1_i) \leq d(L^2_i), \label{dom:3}\\
&\mathcal{V}(L^1_i) \supseteq \mathcal{V}(L^2_i). \label{dom:4}
\end{align}
The dominated labels can be safely discarded during the label extension to speed up the algorithm.

\subsection{Bounded Bidirectional Search}
Bounded bidirectional search partitions the extension of the label-setting algorithm into the \textit{forward extension} and the \textit{backward extension} according to the consumption of a chosen critical resource. In the forward extension, labels are extended in the same way as the original algorithm, while in the backward extension, labels are extended in the reverse direction. After both the forward extension and the backward extension finish, the resulted forward labels and the backward labels join together to generate the complete feasible routes, from which an optimal route can be extracted. Because the number of labels generated usually increases exponentially with the consumption of the critical resource, bounded bidirectional search has the potential to reduce the number of labels generated so as to speed up the algorithm.

\textit{Label Definition.} Let $L^b_i = (\bar{c}(L^b_i), l(L^b_i), n(L^b_i), d(L^b_i), \mathcal{V}(L^b_i))$ be a backward label associated with location $i\in \mathcal{I}^n \cup \{0\}$ which represents a partial path from location $i$ to the depot, where
\begin{itemize}
  \item $\bar{c}(L^b_i)$ is the reduced cost of the path;
  \item $l(L^b_i)$ is the latest arrival time at locations $i$;
  \item $n(L^b_i)$ is the number of locations visited in the path;
  \item $d(L^b_i)$ is the total demand of the visited locations;
  \item $\mathcal{V}(L^b_i)$ is the set of locations that the path can extend to.
\end{itemize}

\textit{Extension functions.} The backward extension starts with an initial label $L^b_0 = \{0, T_{\text{max}}, 0, 0, \mathcal{I}^n \cup \{0\}\}$. For a pair locations $i$ and $j$ and a label $L^b_i$ associated with location $i$, label $L^b_i$ can be extended to location $j$ to create a new label $L^b_j$ by the following extension functions:
\begin{align}
&l(L^b_j) = \min\{L_{\text{max}}, l(L^b_i) - s^n_j - t^n_{ji}\}, \label{bext:1}\\
&n(L^b_j) = n(L^b_i) + 1, \label{bext:2}\\
&\bar{c}(L^b_j) = \bar{c}(L^b_i) + n(L^b_i)(s^n_j + t^n_{ji}) + \bar{c}^n_{ji}, \label{bext:3}\\
&d(L^b_j) = d(L^b_i) + q^n_j, \label{bext:4}\\
&\mathcal{V}(L^b_j) = \left\{k \in \mathcal{V}(L^b_i) \backslash \{j\}~|~d(L^b_j) + q^n_k \leq Q, t^n_{0k} + s^n_k + t^n_{kj} \leq l(L^b_j)\right\}. \label{bext:5}
\end{align}

\textit{Dominance rules.} Let $L^{1b}_i$ and $L^{2b}_i$ be two labels associated with location $i$. Then $L^{1b}_i$ dominates $L^{2b}_i$ if they satisfy the following conditions:
\begin{align}
&\bar{c}(L^{1b}_i) \leq \bar{c}(L^{2b}_i), \label{bdom:1}\\
&l(L^{1b}_i) \geq l(L^{2b}_i), \label{bdom:2}\\
&n(L^{1b}_i) \leq n(L^{2b}_i), \label{bdom:3}\\
&d(L^{1b}_i) \leq d(L^{2b}_i), \label{bdom:4}\\
&\mathcal{V}(L^{1b}_i) \supseteq \mathcal{V}(L^{2b}_i).\label{bdom:5}
\end{align}

\textit{Label Combination.} For a forward label $L^f_i$ and a backward label $L^b_i$ associated with location $i$, they can be joined together to generate a feasible route if
\begin{align}
&\overline{\mathcal{V}}(L^f_i) \cap \overline{\mathcal{V}}(L^b_i) = \{i\}, \label{lc:1}\\
&q(L^f_i) + q(L^b_i) \leq Q + q^n_i, \label{lc:2}\\
&e(L^f_i) \leq l(L^b_i), \label{lc:3}
\end{align}
where $\overline{\mathcal{V}}(L^f_i)$ and $\overline{\mathcal{V}}(L^b_i)$ be the set of locations visited by labels $L^f_i$ and $L^b_i$, respectively. Conditions (\ref{lc:1}), (\ref{lc:2}) and (\ref{lc:3}) guarantee the satisfaction of the elementary constraint, the capacity constraint and the duration constraint, respectively. The reduced cost of the resulted route is equal to  $\bar{c}(L^f_i) + \bar{c}(L^b_i) + e(L^f_i) n(L^b_i)$.

\textit{Critical resource.} The duration $T_{\text{max}}$ is chosen as the critical resource. The break point of the forward extension and the backward extension is determined dynamically as follow. Let $T^f$ and $T^b$ be the upper limits of the critical resource consumption for the forward extension and the backward extension, respectively. That is, for a forward label $L^f$, if $e(L^f) > T^f$, it will not be extended anymore. Similarly, a backward label $L^b$ will not be extended if $l(L^b) < T^b$. Let $\Delta = T_{\text{max}}/16$. Initially, $T^f = \Delta$ and $T^b = T_{\text{max}} - \Delta$. The forward extension and the backward extension are executed. If the number of forward labels is greater than that of the backward extension, $T^f$ remains unchanged and $T^b$ is updated to $T^b - \Delta$. Otherwise, $T^f$ is updated to $T^f + \Delta$ and $T^b$ remains unchanged. This process repeats until $T^f > T^b$.

\subsection{Ng-route Relaxation}
The ng-routes are non-elementary routes first introduced by \citet{baldacci2011new}. To each location $i \in \mathcal{I}^n$, we associate a set of locations $\widehat{\mathcal{I}}^n_i$ $(i \in \widehat{\mathcal{I}}^n_i)$ called {\em ng-set}, e.g. including the closest locations to location $i$. Consider a path $p = (v_1,v_2,\ldots,v_m)$, if $v_1 = v_m$ and $v_m \in \bigcap_{k = 1,\ldots,m - 1}\widehat{\mathcal{I}}^n_{v_k}$, the path $p$ is referred to as an {\em ng-cycle}. A route is an ng-route if it does not contain any ng-cycles. The complexity of the pricing problem based on the ng-route relaxation and the quality of the obtained primal lower bounds depend on the size of the ng-sets. The larger the size of the ng-sets, the better are the lower bounds, but the more difficult it is to solve the pricing problem. Note that if $\widehat{\mathcal{I}}^n_{i} = \mathcal{I}^n~\forall~ i \in \widehat{\mathcal{I}}^{n}$, the pricing problem is equivalent to the ESPPRC, and if $\widehat{\mathcal{I}}^n_{i} = \{i\}$, the pricing problem is equivalent to the non-elementary SPPRC. \citet{martinelli2014efficient} propose an iterative approach to speed up label-setting algorithms based on the ng-route relaxation. First, the elementary constraint is relaxed, and each location is initialized with a small ng-set. Then the label-setting algorithm is called to determine the optimal ng-route. If the optimal ng-route does not contain any ng-cycles, it is elementary and optimal for the original pricing problem. Otherwise the ng-sets of the nodes in the ng-cycles are enlarged to forbid the ng-cycles, and the label-setting algorithm is called again. This process is repeated until an elementary optimal route is found. To handle the ng-route relaxation, the extension functions (\ref{ext:4}) and (\ref{bext:5}) are modified as follows:
\begin{align}
&\mathcal{V}(L^f_j) = (\mathcal{V}(L^f_i) \cup \widetilde{\mathcal{I}}^n_j) \cap \left\{k \in \mathcal{I}^n \backslash \{j\}~|~d(L^f_j) + q^n_k \leq Q,e(L_j) + s^n_j + t^n_{jk} \leq \min\{L_{\text{max}}, T_{\text{max}} - s^n_k - t^n_{k0}\}\right\},\\
&\mathcal{V}(L^b_j) = (\mathcal{V}(L^b_i) \cup \widetilde{\mathcal{I}}^n_{j}) \cap \left\{k \in \mathcal{I}^n \backslash\{j\}~|~d(L^b_j) + q^n_k \leq Q, t^n_{0k} + s^n_k + t^n_{kj} \leq l(L^b_j)\right\},
\end{align}
where $\widetilde{\mathcal{I}}^n_{j} = \mathcal{I}^n \backslash \widehat{\mathcal{I}}^n_{j}$.

\subsection{Label Pruning}

Label pruning is conducted during the extension phase in attempt to identify and drop the labels which  cannot be extended to complete routes with negative reduced costs. Let $\hat{c}^f(i,t)$ be the minimum reduced cost of the partial routes from the depot to location $i \in \mathcal{I}^n$ with arrival time no later than $t$. Then $\hat{c}^f(i,t)$ can be computed by the following dynamic programming:
\begin{align}
&\hat{c}^f(i,t) =
\begin{cases}
0, &\quad \text{if}~ i = 0, 0 \leq t \leq T_{\text{max}}\\
\min_{j \in \mathcal{I}^n \cup\{0\} \backslash\{i\}, w + s^n_j + t^n_{ji} \leq t}\hat{c}^f(j,w) + t + \bar{c}^n_{ji}, &\quad \text{if}~i \neq 0, 0 \leq t \leq T_{\text{max}}.
\end{cases}
\end{align}
Let $\hat{c}^b(i,m)$ be the minimum reduced cost of the partial routes from location $i \in \mathcal{I}^n$ to the depot with the number of the visited locations equal to $m$. Then $\hat{c}^b(i,m)$ can be computed by the following dynamic programming:
\begin{align}
&\hat{c}^b(i,m) =
\begin{cases}
0, &\quad \text{if}~ i = 0, m = 0\\
\infty, &\quad \text{if}~ i = 0, m = 1,\ldots, |\mathcal{I}^n|\\
\min_{j \in \mathcal{I}^n \cup \{0\}\backslash\{i\}}\hat{c}^b(j, m - 1) + (m - 1)t^n_{ji} + \bar{c}^n_{ji}, &\quad \text{if}~i \neq 0, m = 1,\ldots,|\mathcal{I}^n|.
\end{cases}
\label{bprune:1}
\end{align}
The label pruning is done based on $\hat{c}^f(i,t)$ and $\hat{c}^b(i,m)$. A forward $L^f_i$ label can be pruned if
\begin{align}
\bar{c}(L^f_i) + \min_{m + n(L^f_i) \leq |\mathcal{I}^n|}\left\{\hat{c}^b(i, m) + e(L^f_i)m \right\}> 0, \label{bprune:2}
\end{align}
where $n(L^f_i)$ is the number of locations visited by $L^f_i$. Similarly, a backward label $L^b_i$ can be pruned if
\begin{align}
\bar{c}(L^b_i) + \min_{t \leq l(L^b_i)}\left\{\hat{c}^f(i, t) + n(L^b_i)t \right\}> 0.
\end{align}

\subsection{Label-Setting Algorithm for the Part-Time Drivers}
The label-setting algorithm for the full-time drivers can be applied to the part-time drivers with following slight modifications.
\begin{itemize}
    \item The trip duration $T_{\text{max}}$ is set to $L_{\text{max}} + \max_{i \in \mathcal{I}^n}s^n_i + t^n_{i0}$.
    \item The cost $\bar{c}^n_{ij}$ of locations pair $i$ and $j$ is set as follows:
    \begin{align}
&\bar{c}^n_{ij} =
\begin{cases}
-\frac{\nu_i}{2} - \frac{\nu_j}{2}, &\quad \text{if}~i,j \neq 0\\
- \frac{\nu_j}{2}, &\quad \text{if}~i = 0, j \neq 0\\
-\frac{\nu_i}{2} , &\quad \text{if}~ i \neq 0, j = 0.
\end{cases}
\end{align}
    \item Extension function (\ref{ext:2}) is replaced by the following function:
    \begin{align}
    &\bar{c}(L_j) =
\begin{cases}
\bar{c}(L_i) + e(L_j) + \bar{c}^n_{ij}, &\quad \text{if}~j \neq 0\\
\bar{c}(L_i) + \rho e(L_j) + \bar{c}^n_{ij}, &\quad \text{if}~j = 0.
\end{cases}
    \end{align}
    \item Initial backward extension label $L^b_0$ is initialized as $\{0, T_{\text{max}}, \rho, 0, \mathcal{I}^n \cup \{0\}\}$.
    \item Dynamic programming (\ref{bprune:1}) is modified as follows:
    \begin{align}
&\hat{c}^b(i,m) =
\begin{cases}
0, &\quad \text{if}~ i = 0, m = \rho\\
\infty, &\quad \text{if}~ i = 0, m = \{0,\ldots, |\mathcal{I}^n| + \rho\} \backslash \{\rho\}\\
\min_{j \in \mathcal{I}^n \cup \{0\}\backslash\{i\}}\hat{c}^b(j, m - 1) + (m - 1)t^n_{ji} + \bar{c}^n_{ji}, &\quad \text{if}~i \neq 0, m = 1,\ldots,|\mathcal{I}^n| + \rho.
\end{cases}.
\end{align}
    \item Condition (\ref{bprune:2}) is modified as follows:
    \begin{align}
\bar{c}(L^f_i) + \min_{m + n(L^f_i) \leq |\mathcal{I}^n| + \rho}\left\{\hat{c}^b(i, m) + e(L^f_i)m \right\}> 0.
\end{align}
\end{itemize}

\section{Route Enumeration Algorithm}
\label{sec:reap}

\begin{proof}{Proof of Lemma \ref{lemma:re}}
For a route $r \in \mathcal{R}_n$, let $\textbf{SF1}_r$ denote problem \textbf{SF1} where the constraints and the set of feasible routes are restricted to the nodes in $\{i \in \mathcal{I}^n~|~ \alpha_{i,r} = 0\}$, and $\phi(\textbf{SF1}_r)$ be the optimal cost of problem $\textbf{SF1}_r$. Given an upper bound $ub$ of problem \textbf{SF1} and a route $r \in \mathcal{R}_n$, $r$ can not be in any optimal solutions of problem \textbf{SF1} if it satisfies the following condition:
\begin{align}
&c_r + \phi(\textbf{SF1}_r) > ub. \label{con:1}
\end{align}
It is costly to check condition (\ref{con:1}) because a VRP has to be solved for each route. Therefore, we use a lower bound of problem $\textbf{SF1}_r$ instead. Let $\mathcal{I}^n_{r} = \{i \in \mathcal{I}^n~|~ \alpha_{i,r} = 0\}$ and $\mathcal{R}^{n}_r \subseteq \mathcal{R}^n$ be the set of feasible routes which cover locations only in $\mathcal{I}^{n}_r$. A valid lower bound is the LP relaxation of problem $\textbf{SF1}_r$ defined as:
\begin{align}
\textbf{RF1}_r:  \quad \min &  \sum_{r' \in \mathcal{R}^{n}_r}c_{r'}\theta_{r'} \label{rf1r:obj}\\
s.t. \ &\sum_{r' \in \mathcal{R}^{n}_r}\beta_{n',r'}\theta_{r'} \leq \sum_{k = 0}^{\bar{K}}kz^{n'}_k - \beta_{n',r},\quad \forall n'=n+1,\dots, N, \label{rf1r:con1}\\
&\sum_{r' \in \mathcal{R}^{n}_r}\theta_{r'} \leq \bar{K} - \zeta^n_n - 1, \label{rf1r:con2}\\
&\sum_{r' \in \mathcal{R}^{n}_r} \alpha_{i,r'}\theta_{r'} = 1, \quad \forall i\in \mathcal{I}^{n}_r, \label{rf1r:con3}\\
&\theta_{r'} \geq 0, \quad \forall r' \in \mathcal{R}^{n}_r \label{rf1r:con4}
\end{align}

Let $(\pmb{\bar{\mu}}, \pmb{\bar{\nu}})$ be an dual optimal solution of problem $\textbf{RF1}$. Then $(\pmb{\bar{\mu}}, \{\bar{\nu}_i\}_{i \in \mathcal{I}^n_{r}})$ is a dual feasible solution of problem $\textbf{RF1}_r$, and hence $\sum_{n' = n + 1}^{N}\sum_{k = 0}^{\bar{K}}\bar{\mu}_{n'}k\bar{z}^{n'}_k - \sum_{n' = n + 1}^{N}\beta_{n',r}\bar{\mu}_{n'} + \bar{\mu}_n(\bar{K} - \zeta^n_n - 1) + \sum_{i \in \mathcal{I}^n_{r}}\bar{\nu}_i$ is valid lower bound of $\phi(\textbf{SF1}_r)$. If we replace this lower bound in condition (\ref{con:1}), we have
\begin{align}
&c_r + \sum_{n' = n + 1}^{N}\sum_{k = 0}^{\bar{K}}\bar{\mu}_{n'}k\bar{z}^{n'}_k - \sum_{n' = n + 1}^{N}\beta_{n',r}\bar{\mu}_{n'} + \bar{\mu}_n(\bar{K} - \zeta^n_n - 1) + \sum_{i \in \mathcal{I}^n_{r}}\bar{\nu}_i > ub \\
\Leftrightarrow~&c_r - \sum_{i \in \mathcal{I}^n}\alpha_{i,r}\bar{\nu}_i -  \sum_{n' = n + 1}^{N}\beta_{n',r}\bar{\mu}_{n'} - \bar{\mu}_n > ub - \phi(\textbf{RF1}).
\end{align}
\halmos
\end{proof}

The label-setting algorithm introduced in Section \ref{sec:ls} is used to enumerate the target routes by replacing the dominance rules (\ref{dom:4}) and (\ref{bdom:5}) by the following dominance rules:
\begin{align}
&\overline{\mathcal{V}}(L^1_i) = \overline{\mathcal{V}}(L^2_i), \label{re:1}\\
&\overline{\mathcal{V}}(L^{1b}_i) = \overline{\mathcal{V}}(L^{2b}_i). \label{re:2}
\end{align}
The intuition of dominance rules (\ref{re:1}) and (\ref{re:2}) are that two labels cannot dominate each other if the corresponding paths of the labels do not visit the same locations. This ensures all feasible routes can be enumerated. Meanwhile, the ng-route relaxation is not used to speed up the label-setting algorithm. The detailed pseudocode of the route enumeration algorithm is summarized in Algorithm \ref{alg:ea}.

\begin{algorithm}[!h]
	\caption{{\em Route Enumeration}} \label{alg:ea}
	\begin{footnotesize}
		\begin{algorithmic}[1]
			\STATE $ub \leftarrow \infty$, $\delta \leftarrow StepSize$;
			\STATE Solve problem \textbf{RF1} and obtain the optimal primal and dual solutions $\pmb{\bar{\theta}}$ and $(\pmb{\bar{\mu}}, \pmb{\bar{\nu}})$;
			\WHILE{Time limit has not reached \textit{TimeLimit}}
			\STATE $\delta \leftarrow \min\{\delta, ub - \phi(\textbf{RF1})\}$;
			\STATE Enumerate the feasible routes with reduced costs no larger than $\delta$ for dual variables $(\pmb{\bar{\mu}},\pmb{\bar{\nu}})$; \label{line:9}
			\STATE Initialize problem \textbf{SF1} by the enumerated routes and solve it by an MIP Solver;
			\IF{there exists an optimal solution of \textbf{SF1}, i.e., $\pmb{\hat{\theta}}$ and $\phi(\textbf{SF1}) < ub$}
			\STATE $ub \leftarrow \phi(\textbf{SF1})$;
			\IF{$\delta \geq ub - \phi(\textbf{RF1})$} \label{line:14}
			\STATE Break;
			\ENDIF
			\ELSE
			\STATE $\delta \leftarrow \delta + StepSize$;
			\ENDIF
			\ENDWHILE
		\end{algorithmic}
	\end{footnotesize}
\end{algorithm}

\section{Temporal Demand Distribution for Synthetic Instances} \label{sec:distribution}
The following figures present the number of potential customer locations $I_n$ as a function of period $n=1,\dots, N$ for $N=10$ and $N=20$, which controls the temporal variability in generating the synthetic instances.
\begin{figure}[htbp]
	\subfloat[$N=10$\label{sfig:uniformt3}]{
		\includegraphics[width=0.45\linewidth]{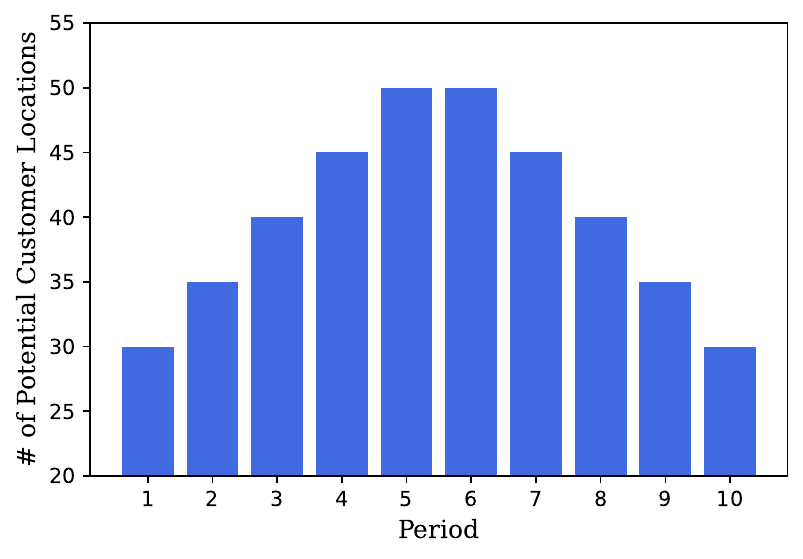}	}\hfill
	\subfloat[$N=20$\label{sfig:uniformwt3}]{
		\includegraphics[width=0.45\linewidth]{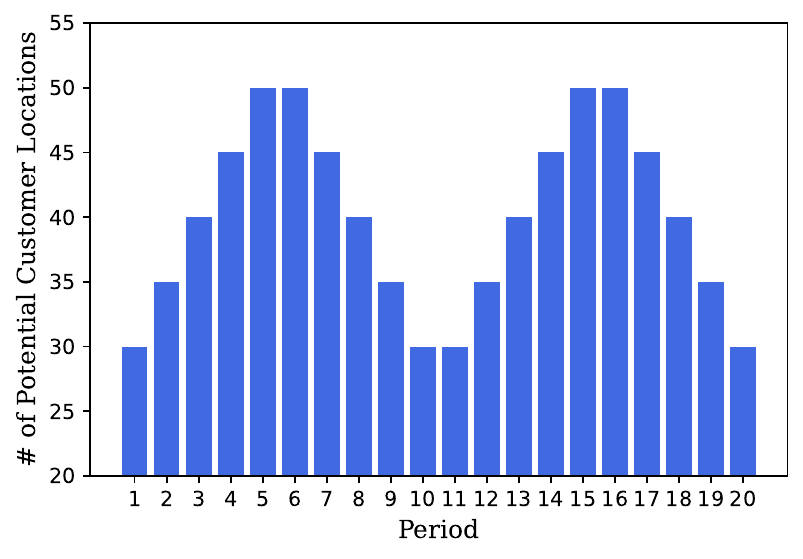}
	}\hfill
	\caption{Average Number of Potential Customer Locations in Different Periods}\label{fig:distribution}
\end{figure}

\section{Additional Numerical Results for AJRP} \label{sec:mintravel}
We compare AJRP to two additional heuristic policies based on driver travel time minimization: (i) MinTravel\_1: the dispatching and routing decisions are obtained by minimizing the total driver time (including return trip to the depot) while satisfying the hard deadline constraints for each epoch; (ii): MinTravel\_2: this is similar to MinTravel\_2 except that the return trip time is excluded in the objective function.  These two heuristics, although myopic, tend to reserve vehicles for future deliveries. Figure \ref{fig:dispatch_mintravel} presents the performance comparison results of AJRP versus these two heuristics on the real data set. We observe AJRP outperforms MinTravel\_1 and MinTravel\_2 consistently by a large margin. This performance gap shows that minimizing driver time is not fully aligned with the system objective: to achieve faster deliveries for realized orders while maintaining a sufficient capacity for future deliveries.  Specifically, these two heuristics dispatch too few drivers even when the capacity is not severely constrained. 
\begin{figure}[htbp]
	\centering
	\includegraphics[width=0.75\linewidth]{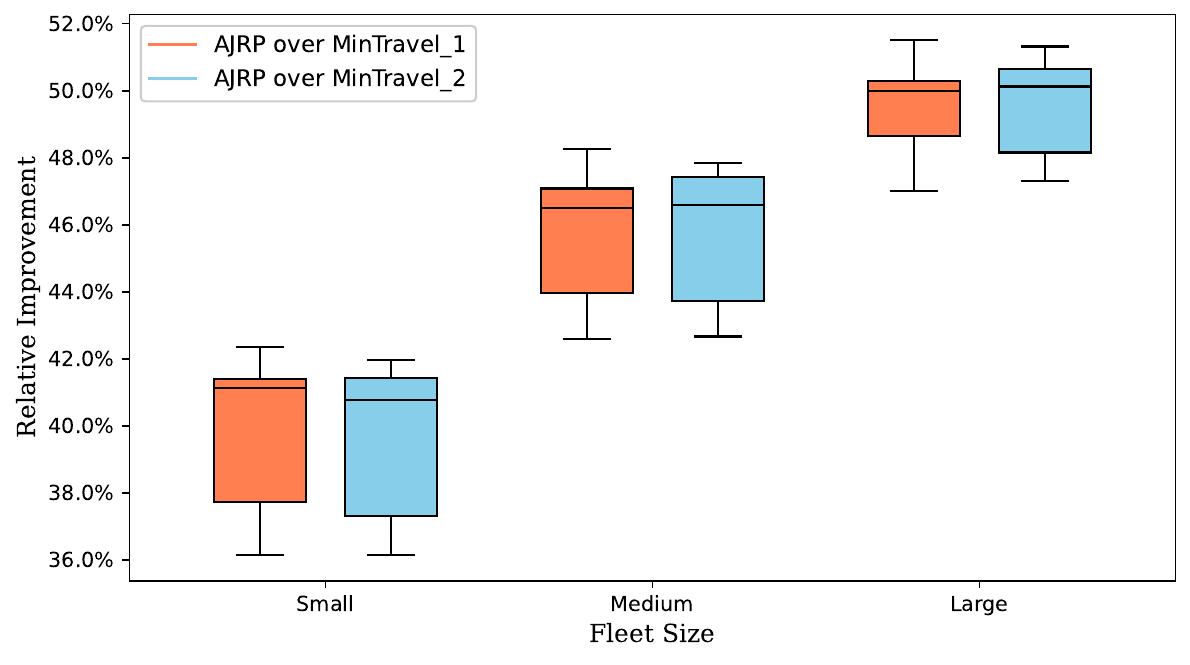}\hfill
	\caption{Performance Evaluation of AJRP, MinTravel\_1, and MinTravel\_2 }\label{fig:dispatch_mintravel}
\end{figure}


\section{Detailed Computational Results on the MTRP Instances} \label{sec:mtrp}
This section presents the detailed computational results of our route enumeration algorithm on public MTRP instances against the state-of-the-art benchmarks. We observe that the proposed algorithm reduces the solution time substantially,  yielding more than 60\% improvement over the benchmarks in many cases. In addition, we  implement a branch-and-price algorithm to solve problem SF1. The branch-and-price algorithm uses two types of branching rules in a hierarchy. First, if the number of routes in the optimal solution is fractional, the algorithm branches on the number of routes. Otherwise, the algorithm selects the arc with the most fractional value and branches on the selected arc.

\begin{table}[!htp]
	\renewcommand{\arraystretch}{1.2}
\caption{Computational Results on the LQL Instances.}
\label{table:set1}
\begin{threeparttable}[t]
\begin{center}
\scalebox{0.75}{
    \begin{tabular}{cccccccccc}
    \toprule
    \multirow{2}[0]{*}{Name} & \multirow{2}[0]{*}{ Number of Customers} &       & \multicolumn{3}{c}{Route Enumeration Algorithm} &       & \multicolumn{1}{c}{\citet{nucamendi2016mixed}} & & \multicolumn{1}{c}{\citet{muritiba2021branch}} \\
    \cline{4-6} \cline{8-8} \cline{10-10}
          &       &       &  Cost & Time & Route Number &       & Time  & & Time \\
    \midrule
    brd14051 & 30    &       & 85729.72 & 0.88  & 746   &       & 4.66  &       & 9.41 \\
    d15112 & 30    &       & 229130.8 & 0.86  & 682   &       & 9.58  &       & 3.67 \\
    d18512 & 30    &       & 84818.95 & 0.77  & 584   &       & 6.85  &       & 2.64 \\
    fnl4461 & 30    &       & 48287.5 & 0.43  & 760   &       & 2.83  &       & 3.65 \\
    nrw1379 & 30    &       & 29994.97 & 0.80  & 795   &       & 5.83  &       & 9.9 \\
    pr1002 & 30    &       & 164763.9 & 0.62  & 1017  &       & 8.09  &       & 4.76 \\
    brd14051 & 40    &       & 106145.5 & 1.16  & 1130  &       & 30.37 &       & 43.49 \\
    d15112 & 40    &       & 290748.6 & 1.78  & 1025  &       & 24.96 &       & 4.88 \\
    d18512 & 40    &       & 107813.4 & 0.89  & 1322  &       & 11.74 &       & 8.21 \\
    fnl4461 & 40    &       & 62672.49 & 0.76  & 967   &       & 17.21 &       & 21.92 \\
    nrw1379 & 40    &       & 36379.77 & 1.04  & 1775  &       & 20.66 &       & 12.54 \\
    pr1002 & 40    &       & 217071.6 & 1.15  & 1615  &       & 52.12 &       & 41.26 \\
    brd14051 & 50    &       & 127480.6 & 1.85  & 1660  &       & 51.54 &       & 73.29 \\
    d15112 & 50    &       & 356433.9 & 3.29  & 2715  &       & 53.75 &       & 33.37 \\
    d18512 & 50    &       & 132115.4 & 1.65  & 1466  &       & 75.64 &       & 49.13 \\
    fnl4461 & 50    &       & 76642.75 & 1.28  & 1216  &       & 32.11 &       & 20.14 \\
    nrw1379 & 50    &       & 45169.23 & 1.89  & 5563  &       & 79.66 &       & 189.44 \\
    pr1002 & 50    &       & 268605.3 & 2.61  & 1651  &       & 60.7  &       & 36.73 \\
    \bottomrule
    \end{tabular}%
  \label{tab:computational_lql}%
}
\end{center}
\end{threeparttable}%
\end{table}

\begin{table}[!htp]
\renewcommand{\arraystretch}{1.2}
\caption{Computational Results on the E Instances.}
\label{table:set2}
\begin{threeparttable}[t]
\begin{center}
\scalebox{0.85}{
    \begin{tabular}{ccccccccc}
    \toprule
    \multirow{2}[0]{*}{Name} &       & \multicolumn{3}{c}{Route Enumeration Algorithm} &       & \multicolumn{1}{c}{\citet{nucamendi2016mixed}} & & \multicolumn{1}{c}{\citet{muritiba2021branch}} \\
    \cline{3-5} \cline{7-7} \cline{9-9}
          &       &  Cost & Time & Route Number &       & Time  & & Time \\
    \midrule
    E-n22-k4 &       & 819.39  & 1.42  & 106   &       & 4.90  &       & 4.32  \\
    E-n23-k3 &       & 1555.87  & 3.01  & 1695  &       & 9.23  &       & 1.37  \\
    E-n30-k3 &       & 1871.08  & 26.63  & 6621  &       & 119.11  &       & 1100.84  \\
    E-n30-k4 &       & 1643.30  & 2.80  & 1023  &       & 22.95  &       & 301.07  \\
    E-n33-k4 &       & 2819.43  & 5.43  & 1890  &       & 24.29  &       & 414.74  \\
    E-n51-k5 &       & 2209.64  & 29.39  & 6127  &       & 2347.51  &       &  t. lim. \\
    E-n76-k7 &       & 2945.25  & 1075.67  & 545451  &       & -     &       &  t. lim. \\
    E-n76-k8 &       & 2677.39  & 152.05  & 81838  &       & -     &       &  t. lim. \\
    E-n76-k10 &       & 2310.09  & 38.19  & 5964  &       & 1700.64  &       & 2797.17  \\
    E-n76-k14 &       & 2005.40  & 19.59  & 5816  &       & 236.64  &       & 271.04  \\
    E-n76-k15 &       & 1962.47  & 18.40  & 5573  &       & 105.43  &       & 109.30  \\
    E-n101-k14 &       & 2922.82  & 152.63  & 36952  &       & -     &       & t. lim. \\
    \bottomrule
    \end{tabular}%
  \label{tab:computational_e}%
}
\end{center}
\begin{tiny}
     \begin{tablenotes}
     \item - means the corresponding instance was not tested.
     \item t. lim. means the computational time reached the time limit.
   \end{tablenotes}
\end{tiny}
\end{threeparttable}%
\end{table}

\begin{table}[!htp]
	\renewcommand{\arraystretch}{1.2}
\caption{Computational Results on the P Instances.}
\label{table:set3}
\begin{threeparttable}[t]
\begin{center}
\scalebox{0.85}{
    \begin{tabular}{ccccccccc}
    \toprule
    \multirow{2}[0]{*}{Name} &       & \multicolumn{3}{c}{Route Enumeration Algorithm} &       & \multicolumn{1}{c}{\citet{nucamendi2016mixed}} & & \multicolumn{1}{c}{\citet{muritiba2021branch}} \\
    \cline{3-5} \cline{7-7} \cline{9-9}
          &       &  Cost & Time & Route Number &       & Time  & & Time \\
    \midrule
    P-n16-k8 &       & 382.90  & 0.46  & 803   &       & 1.90  &       & 0.36  \\
    P-n19-k2 &       & 812.15  & 2.20  & 1711  &       & 9.13  &       & 1.27  \\
    P-n20-k2 &       & 905.19  & 2.29  & 1759  &       & 11.60  &       & 2.07  \\
    P-n21-k2 &       & 937.10  & 4.23  & 1663  &       & 11.07  &       & 6.36  \\
    P-n22-k2 &       & 993.10  & 3.69  & 2067  &       & 11.00  &       & 1.87  \\
    P-n22-k8 &       & 623.40  & 0.72  & 1147  &       & 3.08  &       & 0.29  \\
    P-n23-k8 &       & 561.33  & 0.78  & 1328  &       & 2.68  &       & 0.29  \\
    P-n40-k5 &       & 1537.79  & 9.30  & 1183  &       & 213.28  &       & 196.62  \\
    P-n45-k5 &       & 1912.31  & 18.12  & 2109  &       & 495.82  &       & 3179.19  \\
    P-n50-k7 &       & 1547.89  & 10.45  & 1766  &       & 117.48  &       & 211.15  \\
    P-n50-k8 &       & 1448.92  & 8.62  & 1534  &       & 185.76  &       & 958.75  \\
    P-n50-k10 &       & 1296.48  & 5.86  & 1724  &       & 112.32  &       & 13.57  \\
    P-n51-k10 &       & 1419.43  & 6.32  & 2062  &       & 84.17  &       & 82.91  \\
    P-n55-k7 &       & 1766.56  & 16.60  & 2068  &       & 790.31  &       & 3033.76  \\
    P-n55-k8 &       & 1614.61  & 13.99  & 8909  &       & 170.89  &       & 136.19  \\
    P-n55-k10 &       & 1438.60  & 7.86  & 5961  &       & 117.69  &       & 15.53  \\
    P-n55-k15 &       & 1280.92  & 5.83  & 1678  &       & 48.01  &       & 4.23  \\
    P-n60-k10 &       & 1676.35  & 12.45  & 2723  &       & 620.81  &       & 674.62  \\
    P-n60-k15 &       & 1462.50  & 7.35  & 2247  &       & -     &       & 13.14  \\
    P-n65-k10 &       & 1928.46  & 17.91  & 2566  &       & 915.83  &       & 6762.68  \\
    P-n70-k10 &       & 2097.17  & 22.35  & 3942  &       & 1415.78  &       & 5877.93  \\
    P-n76-k4 &       & -     & t. lim. & 297700  &       & t. lim. &       & t. lim. \\
    P-n76-k5 &       & -     & t. lim. & 129279  &       & t. lim. &       & t. lim.  \\
    \bottomrule
    \end{tabular}%
  \label{tab:addlabel}%
}
\end{center}
\begin{tiny}
     \begin{tablenotes}
     \item - means the corresponding instance was not tested.
     \item t. lim. means the computational time reached the time limit.
   \end{tablenotes}
\end{tiny}
\end{threeparttable}%
\end{table}


Table \ref{table:singlePeriod2} compares the computational performance between the route enumeration algorithm and the branch-and-price algorithm. We observe that the route enumeration algorithm outperforms the branch-and-price algorithm consistently. The main reason is that the branch-and-price algorithm has to explore a number of nodes to finally achieve the optimal solution.
 
\begin{table}[!htp]
\renewcommand{\arraystretch}{1.2}
\caption{Computational Comparison of the Route Enumeration and  Branch-and-Price Algorithms on the MTRP Instances}
\label{table:singlePeriod2}
\begin{center}
\scalebox{0.85}{
    \begin{tabular}{cccccccc}
    \toprule
    \multirow{3}[0]{*}{Class} & \multirow{3}[0]{*}{\parbox{1.5cm}{~~~Total\\ Instances}} & \multicolumn{2}{c}{Route Enumeration Algorithm} & & \multicolumn{3}{c}{Branch-and-Price}\\
    \cline{3-4} \cline{6-8}
          &       & Instances  & Average Time & & Instances  & Average Time & Number of \\
          &       & Tested/Solved & (In seconds) & & Tested/Solved & (In seconds) & nodes \\
    \midrule
    LQL  & 180   & 180/180 & 1.32  & & 180/180 & 22.86  & 17.6  \\
    E  & 12    & 12/12 & 127.10  & & 12/9   & 496.43  & 24.11   \\
    P  & 23    & 23/21 & 8.45  & & 23/21 & 276.08  & 20.80  \\
    \bottomrule
    \end{tabular}%
}
\end{center}
\end{table}

\end{APPENDICES}
\end{document}